\documentclass{amsart}
\usepackage{amssymb}
\usepackage{amsfonts}
\usepackage{amsbsy}

\usepackage{latexsym}
\usepackage[mathscr]{eucal}
\usepackage{verbatim}

\theoremstyle{plain}
\newtheorem{theorem}{Theorem}[section]
\newtheorem{corollary}[theorem]{Corollary}
\newtheorem{lemma}[theorem]{Lemma}
\newtheorem{proposition}[theorem]{Proposition}
\newtheorem{remark}[theorem]{Remark}
\newtheorem{definition}[theorem]{Definition}
\newtheorem{example}[theorem]{Example}
\newtheorem{question}[theorem]{Question}
\def\sideremark#1{\ifvmode\leavevmode\fi\vadjust{\vbox to0pt{\vss
\hbox to 0pt{\hskip\hsize\hskip1em
\vbox{\hsize2cm\tiny\raggedright\pretolerance10000
\noindent#1\hfill}\hss}\vbox to8pt{\vfil}\vss}}}

\newcommand{\be}{\begin{equation}\label}
\newcommand{\ee}{\end{equation}}
\newcommand{\bq}{\begin{equation*}}
\newcommand{\eq}{\end{equation*}}
\newcommand{\ba}{\begin{align*}}
\newcommand{\ea}{\end{align*}}
\newcommand{\bp}{\begin{proof}}
\newcommand{\ep}{\end{proof}}
\newcommand{\bL}{\begin{lemma}\label}
\newcommand{\eL}{\end{lemma}}
\newcommand{\bP}{\begin{proposition}\label}
\newcommand{\eP}{\end{proposition}}
\newcommand{\bC}{\begin{corollary}\label}
\newcommand{\eC}{\end{corollary}}
\newcommand{\bT}{\begin{theorem}\label}
\newcommand{\eT}{\end{theorem}}
\newcommand{\bR}{\begin{remark}\label}
\newcommand{\eR}{\end{remark}}
\newcommand{\bD}{\begin{definition}\label}
\newcommand{\eD}{\end{definition}}

\DeclareMathOperator{\tr}{Tr}
\DeclareMathOperator{\diag}{diag}

\begin{document}
\title{Finite sums of projections in von Neumann  algebras.}
\author{Herbert Halpern}
\address{University of Cincinnati\\
          Department of Mathematics\\
          Cincinnati, OH, 45221-0025\\
          USA}
\email{halperhp@ucmail.uc.edu}
\author{Victor Kaftal}
\address{University of Cincinnati\\
          Department of Mathematics\\
          Cincinnati, OH, 45221-0025\\
          USA}
\email{victor.kaftal@math.uc.edu}
\author{Ping Wong Ng}
\address{University of Louisiana\\
          Department of Mathematics\\
          Lafayette, LA, 70504\\
          USA}
\email{png@louisiana.edu}
\author{Shuang Zhang}
\address{University of Cincinnati\\
          Department of Mathematics\\
          Cincinnati, OH, 45221-0025\\
          USA}
   \email{zhangs@ucmail.uc.edu}
\date{July 27, 2010}

\begin{abstract}
We first prove that in a $\sigma$-finite von Neumann factor $M$, a positive element $a$  with properly infinite
range projection $R_a$ is a linear combination of projections with positive coefficients if and only if  the essential norm $\|a\|_e$ with respect to the closed two-sided ideal $J(M)$ generated by the finite projections of $M$ does not vanish. Then we show that if  $\|a\|_e>1$, then $a$ is a finite sum of  projections. Both these results are extended to general properly infinite von Neumann algebras in terms of central  essential spectra. 

Secondly, we provide a necessary condition for a positive operator $a$ to be a finite sum of projections  in terms of the principal ideals generated by the excess part $a_+:=(a-I)\chi_a(1,\infty)$ and the defect part $a_-:= (I-a)\chi_a(0, 1)$ of $a$; this result appears to be new also for $B(H)$. 

Thirdly, we prove that  in a type II$_1$ factor a sufficient  condition for a positive diagonalizable operators to be a finite sum of projections is that $\tau(a_+)- \tau(a_-)>0$.

\end{abstract}
\title{Finite sums of projections in von Neumann algebras.}
\maketitle
\section{Introduction} \label{S:intro}

The goal of this article  is to study the following two related
problems in the context of general von Neumann algebras.

\begin{enumerate}
\item [(A)] Which positive operators  are  linear combinations  of projections with positive
coefficients (called \textit {positive combinations of
projections})?
\item [(B)] Which positive operators  are \textit {finite sums of
projections}?
\end{enumerate}
Let us first give some historical background about problem
(A). Fillmore in \cite {Fp67}
proved that all operators in $B(H)$ are linear combinations of
projections. This result was  extended  to properly infinite von Neumann algebras by Pearcy and Topping
\cite {PcTd67}, to type II$_1$ factors by Fack and De La Harpe \cite
{FdlH80} and by Pearcy and Topping \cite {PcTd69}, and to von Neumann
algebras with no type I direct summands with infinite
dimensional center by Goldstein and Paszkiewicz \cite {GsPa92}.
It was also extended to several  types of C$^*$ algebras by various authors including Marcoux \cite {Mlw02}, \cite{Mlw08} and Marcoux and Murphy  \cite {MM98}.

The first result specifically on problem (A) was obtained by Fillmore who proved in \cite {Fp69}  that positive \textit{invertible} operators in $B(H)$
(for a separable Hilbert space $H$) are positive combinations of projections.  This result later was extended to all von Neumann algebras with no finite type I direct summands by Bikchentaev \cite [Lemma 5] {Ba05}. In this article we
present an alternative proof for the same result and provide an
estimate of the number of needed projections (Proposition \ref
{P2.2}). 

It was remarked by Fillmore  in \cite{Fp67} that infinite
rank compact operators on a separable Hilbert space  cannot be
positive combinations of projections;  Fong and Murphy  in
\cite[Theorem 11] {FoMj85} showed that these operators are the only
exceptions in $B(H)$.  

We prove that for $\sigma-$finite von Neumann factors that are finite or of type III, all positive operators are positive combinations of projections, while for a type II$_\infty$ factor $M$ the only exceptions are, like in $B(H)$, those operators with infinite range projection that belong to the ideal of relative compact operators (Theorem \ref {T:2.11}, Corollary \ref{C3.5}).  For a  general von Neumann algebra $M$, the ideal of relative compact operators $J(M)$ is the closed two-sided ideal generated by the finite projections of $M$ and was first studied  by Breuer \cite {Bm68}, \cite {Bm69} and Sonis \cite {Sm71}.

For global von Neumann algebras,  we obtain a characterization of positive combinations of projections with properly infinite range projection  in terms of a spectral property of the operator (Theorem \ref {T:2.11}). Its connection with ideals is more technical (Theorem \ref {T3.4}) and is formulated in terms of the notion of central essential spectrum relative to a central ideal. This notion was introduced  by Halpern in a cycle of papers including \cite {Hh72}
and \cite {Hh77}, (see also Stratila and Zsido  \cite {SsZl73}). The relevant facts are  summarized in Section \ref {S:central}.

The key ingredient for our characterization both for the factor and non-factor case is Lemma \ref {L:2.7}, new also for $B(H)$, which proves that the direct sum of an arbitrary positive operator and of a locally invertible operator on a ``large" subspace is a positive combination of projections. A  slightly weaker formulation of this lemma holds also for a large class of C$^*$-algebras and  is used in another paper
\cite{KNZC*10}.

\medskip

Next we give some background on problem (B). Fillmore characterized the positive finite rank operators that are sums of finitely many projections:
\vspace{-0.15cm}
\bT{T1.1}\cite [Theorem 1] {Fp69} 
 Let $a\in B(H)^+$ have finite rank. Then $a$ is a sum of projections if and only if $\tr(a)\ge$ rank$(a)$ and $\tr(a)\in \mathbb N$.
\eT
 As reported in a survey article
\cite [Theorem 4.12] {Wpy94} by Wu quoting unpublished joint work of
Choi and Wu in 1988, they  proved that positive operators with essential norm greater than 1 are \textit{finite sums of projections}. This result is presented together with other related results  in their recent paper \cite[Theorem 2.2] {WpCm09}. The special case of $\alpha I$   with $\alpha
>1$  follows also from a delicate analysis  which characterize the numbers $\alpha
$  for which  $\alpha I$ is a sum of at least $n$ projections (see in \cite {KRS02} and 
several other papers.)

In recent years related topics became quite active again.  In \cite {DFKLOW} Dykema,
Freeman, Kornelson, Larson, Ordower, and Weber, motivated by problems in frame theory, proved by  different techniques that all positive operators   with essential norm greater than 1 are   \textit{strong sums of projections}. Their work inspired other papers, including \cite {AMRS}, \cite {KkLd04}, and  work  by three of the authors of the present article (\cite {KNZW*09}, \cite {KNZC*09}, \cite{KNZC*10}).

 \cite {KNZW*09} obtained a characterization of strong sums of projections  in properly infinite $\sigma$-finite von Neumann factors (thus including $B(H)$)  and for \textit{diagonalizable} operators, also in type II factors. These characterizations were given in terms of the \textit{excess}
and \textit{defect} operators $a_+$ and $a_-$ which will play an  important role also in the present paper (see Definition \ref
{D5.3} and Theorem \ref {T5.4} for a precise statement.)

\cite {KNZC*09} presented a characterization of \textit{strict} sums of projections in the multiplier algebras of $\sigma$-unital  purely infinite simple C$^*$ algebras and \cite {KNZC*10} provided sufficient conditions for operators to be finite sums of projections in the multiplier algebras of $\sigma$-unital  purely infinite simple C$^*$ algebras.

The full problem (B), i.e., the complete characterization of the positive operators that are \textit{finite} sums  of projections  is still open even in the $B(H)$ case, but we present some sufficient and  some necessary conditions.

For  $\sigma-$finite factors  a sufficient condition is the natural $B(H)$ analog:  having the essential norm larger than 1. The essential norm is relative to the compact ideal $J(M)$, namely the norm in the quotient algebra $M/J(M)$ (see Corollary \ref{C4.4}). In the type III case, $J(M)=\{0\}$ and hence the condition is for the operator to have norm   larger than 1.

For global algebras we  replace the essential norm  with the ``central essential norm", again in terms of central essential spectra with respect to a certain central ideal (Theorem \ref{T4.3}.) The reader mainly interested  in factors or, in particular, in $B(H)$, can simply assume that all the ``central" objects mentioned are scalars.

Our approach is based on  Lemma \ref {L4.2} where we show that the direct sum  of an arbitrary positive operator  and a scalar multiple   greater than 1  of a  ``large" projection  is always a finite sum of projections.
As a consequence, we obtain in Theorem \ref {T4.3}  that if  the   ``central essential norm" of an operator $a\in M^+$ is greater than and bounded away from the identity, then $a$ is a finite sum of projections. Restricted to the special  case  $B(H)$, our  approach gives an alternative proof to the result of Choi and Wu \cite{WpCm09}, and hence, of Dykema and all \cite {DFKLOW}.

Then in Theorem \ref {T5.6} we find necessary conditions for a positive operator $a$ to be a finite sum of projections. It is formulated in terms
the two-sided (non-closed) principal ideals  of $M$ generated by the
excess and defect operators of $a$.  This condition is new for
$B(H)$ and generalizes some results in \cite{WpCm09} (see Remark
\ref {R5.11} for more details).

From the initial work of Fillmore, Pearcy, and Topping, the finite case 
has proven to be more delicate than the properly infinite case. In \cite[Theorem 1]{KNZW*09} we considered
positive \textit{diagonalizable} operators in a II$_1$ factor $M$
with trace $\tau$  and proved  that the inequality $\tau(a_+)\ge
\tau(a_-)$  is necessary and sufficient for $a$ to be a strong sum
of projections. Via a different and considerably more complicated analysis, we show in the present article (Theorem \ref {T6.5}) that the
strict inequality $\tau(a_+)> \tau(a_-)$ is sufficient also for $a$
to be a finite sum of projections. However, the  problem
remains open when $a$ is not  diagonalizable or when $\tau(a_+)= \tau(a_-)$.

\

We wish to thank  Pei Yuan Wu, Viacheslav Rabanovych, and Kostyantyn
Yusenko for very useful information and references  in this research
area.

\

Notations: \\
For every $x\in \mathbb R$,  we will denote by $\lceil x \rceil $ the smallest integer $n \ge x$ and by $\lfloor x \rfloor $ the
largest integer $n \le x$, i.e., the \textit{the integer part} of $x$.\\
For every self-adjoint element $a$, $\chi_a$ will denote
the (projection valued) spectral measure of $a$ and $R_a$ the range projection of
$a$\\
We shall use standard notations for  von Neumann
algebras. $M$ will be a von Neumann algebra represented on a Hilbert space $H$. In general, we will not assume that  $M$ is  a
factor nor that it is  $\sigma$-finite (i.e., countably
decomposable). For semifinite algebras, $\tau$ will denote a
faithful, normal, semifinite trace. 
$Z=M\cap M'$ will denote the  center of $M$. If $p\in M$, $c(p)$ will denote the central
support of $p$ and $M_p$  the restriction of $pMp$ to the
subspace $pH$. We will identify  $M_p$  with $pMp$. 
\\
$J(M)$ will denote the closed two-sided ideal generated by the finite projections of $M$. If $M=B(H)$, then $J(M)$ is the ideal $K(H)$ of the compact operators on $H$. If $M$ is type III, then $J(M)=\{0\}$.
The essential norm with respect to the $J(M)$, i.e., the quotient norm in $M/J(M)$ is denoted by  $\|\cdot\|_e$.

\section{Positive combinations of projections} \label{S:pos}

Goldstein and Paskiewicz proved in \cite {GsPa92} that a von Neumann algebra has the property that its elements  are linear combinations of projections in the algebra if and only if the algebra has no  finite type I direct summands with infinite dimensional  center, i.e.,  if it is a  direct sum of a properly infinite algebra, a type II$_1$ algebra, and at most a finite number of matrix algebras.

Using his characterization of sums of two projections \cite
[Corollary of Theorem 2]{Fp69},  Fillmore proved in \cite[Corollary of Theorem 3]{Fp69} that every positive \emph{invertible} operator acting on a  separable Hilbert space can be decomposed into a positive combination of projections. 

It is easy to see that this property does not extend to general von Neumann algebras.

\bL{L2.1}
Let $\mathscr A$ be an infinite dimensional abelian von Neumann algebra. Then there are positive invertible operators in $\mathscr A$ that are not linear combinations of projections.
\eL
\bp
Since $\mathscr A$ is infinite dimensional, it contains an infinite sequence $\{e_k\}$ of mutually orthogonal projections with  $\sum_{k=1}^\infty e_k = I$. Let $\mu_k\in [1,2]$ be a sequence with infinitely many distinct entries. Let $a:=\sum_{k=1}^\infty \mu_ke_k$. Then $a\in \mathscr A$ and $a\ge I$, hence it is a positive invertible operator. We claim that $a$ is not a linear combination of projections in $\mathscr A$.  Assume by contradiction that  $a=\sum_{j=1}^n\lambda_j p_j$ for some $\lambda_j\in \mathbb C$ and projections $p_j\in \mathscr A.$ Since the projections commute, we can subdivide them into mutually orthogonal projections  and thus obtain that $a= \sum_{j=1}^{n'}\lambda'_j p'_j$ with $p'_j $ mutually orthogonal projections in $\mathscr A.$ But then, on the one hand we would obtain that  the spectrum Sp$(a)= \{\lambda'_j\}_1^{n'}$ is finite and on the other hand that Sp$(a)=\overline{\{\mu_k\}_1^\infty}$ is infinite, a contradiction. \ep

Abelian von Neumann algebras are  a special case  of finite type I algebras. In \cite[Lemma 5 (C)] {Ba05},
Bikchentaev proved using \cite[Theorem 3]{GsPa92} that in any von Neumann algebra with no finite type I direct summands,  positive invertible operators are always positive combinations of projections.

In the following proposition we present a  proof of  this
result that provides an explicit estimate  of  the  number of projections  required, which we need in Theorem \ref {T6.5} below.  Our proof is an adaptation of  Fong's  $B(H)$ argument in \cite[Lemma 8]{FoMj85},  based on the notion of the following algebra constants $N_o$ and $V_o$, which he introduced for $B(H)$, but which exist also for many other operator algebras:

\bD{D:2.2}
An operator algebra $M$ has constants $N_o$ and $V_o$ if every selfadjoint operator $a\in M$ can be decomposed into a real linear combination $a= \sum_{j=1}^{N_o}\alpha_jp_j$ of $N_o$ projections $p_j\in M$ such that
$$
  \sum_{j=1}^{N_o}|\alpha_j|\le V_o\|a\|.
$$
\eD
The existence and the estimates of the constants $N_o$ and $V_o$  for von Neumann algebras with no finite type I direct summands with infinite dimensional  center are given in  Theorems 1--3  in \cite{GsPa92}, explicitly for $N_o$, and implicitly (i.e., as a simple consequence of the statements of the theorems) for $V_o$:
\begin{align}  &\text{If $M$ is properly infinite then $N_o=6$ and $V_o=8;$}\notag\\
& \text{if $M$ is of type II$_1$, then $N_o=12$ and $V_o=14;$}\label{e:2.1}\\
&\text{if $M$ is the direct sum of $m$ matrix algebras, then $N_o= m+4$ and $V_o= m+4$.}\notag
\end{align}

\bP{P2.2}
Let $M$ be a von Neumann algebra that has no finite type I direct
summands with infinite dimensional center and hence has  constants $N_o$ and $V_o$. 
Assume that  $a\in M$ is positive and is invertible.  Then $a$  is a positive
 combination  of $N_o+\lceil  V_o \Big(\|a\|\|a^{-1}\|-1\Big)\rceil+1$ projections in $M.$
\eP
\bp
Set $\nu:=\frac{1}{\|a^{-1}\|}$. Then $a\ge \nu I$ and $\chi_a[\nu, \|a\|]$=1. Partition the interval $[\nu, \|a\|]$ into $n= \lceil V_o\Big(\frac{  \|a\|}{\nu}-1\Big)\rceil $ equal subintervals with endpoints
 $\lambda_k$  and set  $$e_k:=\begin{cases}\chi_a[\nu, \lambda_1]&k=1\\
 \chi_a(\lambda_{k-1},\lambda_k]&2\le k\le n.\end{cases}\qquad\text{and}\qquad b:=a- \sum_{j=1}^n \lambda_{j} e_j.$$ Then  $\sum_{j=1}^n e_j =I$,  $b\ge 0$ and
$$\|b\|\le \frac{\|a\|-\nu}{n}\le \frac{\nu}{V_o}.$$ By the hypothesis that $M$ has  constants $N_o$ and $V_o$, we can decompose $b= \sum_{i=1}^{N_0}\alpha_i q_i$
into a real linear combination of $N_o$ projections $q_i \in M$ 
satisfying the condition that $ \sum_{i=1}^{N_0}|\alpha_i |\le V_o\|b\|\le \nu$. But then
\ba
a&=  \sum_{j=1}^n (\lambda_j-\nu) e_j + \nu I + \sum_{i=1}^{N_0}\alpha_i q_i\\
&= \sum_{j=1}^n (\lambda_j-\nu) e_j +  \sum _{ \alpha_i\ge 0
}\alpha_i q_i  + \sum_{\alpha_i< 0  }(-\alpha_i)(I- q_i)  +
 \big(\nu - \sum_{\alpha_i<
0  }(-\alpha_i) \big) I.
\end{align*} As desired, the number of projections in the above linear
combination does not exceed  $N_o+n+ 1$ and all the coefficients are
nonegative.
\ep

A positive element $a$ such that  $a\ge \nu R_a$ for some
$\nu>0$ is  invertible in the reduced algebra $M_{R_a}$. We call such an element \textit{locally invertible}.

If $p$ is a projection in $M$, we identify $M_p$  with the hereditary subalgebra $pMp$ of $M$.   Clearly, if $a\in M_p$ is a linear combination of projections in $M_p$, then $a$ is  also a linear combination of projections in $M$. The converse, however, does not hold: $a$ can be a linear combination of projections in $M$,  but  as the restriction to $pH$ of a projection $q$ fails to be a projection unless $q$ commutes with $p$,  it does not follow that $a$ is a linear combination of projections in $M_p$. If however $a=\sum_{j=1}^n\lambda_jp_j$ is a positive combination of projections in $M$, then  $p_j\le \frac{1}{\lambda_j}a\le  \frac{\|a\|}{\lambda_j}p$ and hence $p_j\in M_p$ for all $j$. Thus $a$ is also a positive combination of projections in $M_p$. To summarize:

\bL{L2.3}Let $a\in M^+$. Then $a$ is a positive combination of projections in $M$ if and only if  $a$  is a positive combination of projections in $M_{R_a}$.
\eL

Thus an immediate consequence of Proposition \ref {P2.2} and its proof is the following.

\bC{C2.4} Let $a\in M^+$ and assume that is locally invertible, i.e., $a\ge \nu R_a$ for some $\nu \ge 0$, and that $M_{R_a}$ has constants $N_o$ and $V_o$. Then $a$ is a positive combination of $N_o+ \lceil V_o\big(\frac{\|a\|}{\nu}-1\big)\rceil+1$ projections in $M$.
\eC

Recall that $M_{R_a}$ has constants $N_o$ and $V_o$ if and only if $M_{R_a}$ has no finite type I direct summands with infinite dimensional center. To guarantee that this holds for every $a\in M$, we must assume that $M$ itself has no (finite or infinite) type I direct summands with infinite dimensional  center. 

\bC{C2.5} A von Neumann algebra $M$ has the property that every  positive locally invertible operator in $M$ can be decomposed into a positive combination of projections in $M$ if and only if 
it has no type I direct summands with infinite dimensional  center.
\eC
\bp
Assume that $M$ has no type I direct summands with infinite dimensional center and let $a\in M^+$ be locally invertible. Then $a$ is invertible in the reduced algebra $M_{R_a}$. Now $M_{R_a}$  has also no type I direct summands with infinite dimensional center and hence a fortiori, has no \textit{finite} type I direct summands with infinite dimensional center. Thus the conclusion follows from Proposition \ref
{P2.2}  and Lemma \ref {L2.3}.

For the opposite implication, assume by contradiction that  $M$
has a type I direct summand $M_g$  with infinite dimensional center  $(M\cap M')_g$, for some central projection $g$. By definition  there is an abelian projection $p\in M$ with central support $c(p)=g$ and hence $M_p= (M\cap M')_p$ is an infinite dimensional abelian algebra. By Lemma \ref {L2.1}, $M_p$ contains positive invertible operators that are not positive linear combination of projections in $M_p$. Seen as elements of $M$, they are positive  locally invertible and by Lemma   \ref {L2.3}, they are not positive combinations of projections in $M$.
\ep

\begin{question}\label{Q:2.6}
In the case of  finite type I algebras with infinite dimensional center  we are not aware of conditions that determine which elements   are linear combinations  or are positive combinations  of
projections.
\end{question}

Invertibility is of course not a necessary condition for a positive operator in $B(H)$ to be  a positive combination
of projections -- indeed by the  spectral theorem, all positive finite rank operators are positive combinations of
projections. Fillmore noted that compact operators with infinite rank cannot be a positive combination of projections
 \cite [Remark 5] {Fp67}. Fong and Murphy  showed in  \cite [Lemma 9]{FoMj85} that these are the only exceptions.
 To extend this result to von Neumann algebras, we will will need to follow a different approach. The key steps
 is provided by the following lemma which shows that a direct sum of a positive operator with a positive invertible
 operator of ``sufficiently large"  range is a positive combination of projections.

 \bL {L:2.7}
 Let $e$ and $f$ be orthogonal projections in $M$ with $e\prec f$ and such that $M_{f}$ has no  finite type I direct
  summands with infinite dimensional center.  Let $b=be=eb, d=df=fd$ be positive elements in $M$ such that  $d\ge \nu f$
  for some $\nu>0$. Then $a:=b+d$ is a positive combination of projections in $M.$
 \eL
 \bp
 If $b=0$, then $a=d$ is an invertible operator in the algebra $M_f$, thus  $a$ is positive combination of projections in $M$ by Proposition \ref {P2.2} and Lemma \ref {L2.3}. Assume from now on that $b\not= 0$. 
 
We first prove the statement under the additional assumption that $\nu> \|b\|$.

 Choose a partial isometry $v\in M$ for which $v^*v=e$
and $f':=vv^*\le f.$ With respect to the matrix units $e,v, v^*,
f'$, define the projections
$$ q_-:= \begin{pmatrix}\frac{b}{\|b\|}&-\big(\frac{b}{\|b\|}-\frac{b^2}{\|b\|^2}\big)^{1/2}\\-\big(\frac{b}{\|b\|}
-\frac{b^2}{\|b\|^2}\big)^{1/2}& v(e-\frac{b}{\|b\|})v^*\end{pmatrix}\qquad \text{and}\qquad q_+:= \begin{pmatrix}
\frac{b}{\|b\|}&\big(\frac{b}{\|b\|}-\frac{b^2}{\|b\|^2}\big)^{1/2}\\ \big(\frac{b}{\|b\|}-\frac{b^2}{\|b\|^2}\big)^{1/2}&
 v(e-\frac{b}{\|b\|})v^*\end{pmatrix}.$$
Then $q_-+q_+= \frac{2}{\|b\|}b+ 2f'-\frac{2}{\|b\|}vbv^*$, and
hence
$$
a= \frac{\|b\|}{2}q_-+ \frac{\|b\|}{2}q_+  + d-\|b\|f' + vbv^*.
$$
Since $$d-\|b\|f' + vbv^*\ge d-\|b\|f'\ge(\nu-\|b\|) f,$$
 it follows  by Proposition \ref {P2.2} and Lemma \ref {L2.3}  that  the locally invertible element $d-\|b\|f' + vbv^*$ is a positive combinations of projections
 in $M$, and hence, so is $a$.
 
 Now we remove the assumption that $\nu > \|b\|$. 
First of all, we reduce the proof to the case when $b$ itself is locally invertible. By   \cite [Lemma 3.2] {Hh77}, we can decompose $f=f'+f''$ into the sum of two  projections $f'\sim f''\sim f$ \textit{that commute with $d$.} Then 
$$a= b\chi_b(0, \frac{\nu}{2})+ df'+ b\chi_b[ \frac{\nu}{2}, \infty)+df''.
$$
Now $\chi_b(0, \frac{\nu}{2})\le e\prec f\sim f'$, $\|b\chi_b(0, \frac{\nu}{2})\|\le  \frac{\nu}{2}$ while $df'\ge \nu f'$, and $M_{f'}$, like $M_f$, has no finite type I direct summands with infinite dimensional
center. Thus for the first part of the proof, $ b\chi_b(0, \frac{\nu}{2})+ df'$ is a positive combination of projections. It remains to consider $b\chi_b[ \frac{\nu}{2}, \infty)+df''$, or, to simplify notations, to just assume that $b\ge \frac{\nu}{2} e$.

Now we decompose the central support $c(f)$ of $f$ into the sum of  four central projections $$c(f) = g_1+g_2+g_3+g_4$$(some of which can vanish) as follows\\
$g_1:= c(f)-c(e)$ \quad (notice that  $e\prec f$ implies that  $c(e)\le c(f)$;)\\
$g_2$ is such that $M_{eg_2}$ has no finite type I direct summands with infinite dimensional center and $M_{eg_2^\perp}$
is a finite type I algebra;\\
$g_3f$ is finite;\\
$g_4f$ is properly infinite. \\
Thus $a= \sum_{j=1}^4ag_i$, so it is enough to prove that  $ag_i$ is a positive combination of projections for $i=1,\cdots,4$.

Since $ag_1= d(c(f)-c(e))\ge  \nu f\big(c(f)-c(e)\big)$ and $M_ { f\big(c(f)-c(e)\big)}$ being a direct summand of $M_f$, also has no  finite type I direct summands with infinite dimensional center, the conclusion follows  from Proposition \ref {P2.2} and Lemma \ref {L2.3}.

Since $ag_2= beg_2+ dfg_2$ is a sum of two positive locally invertible operators  and both $M_{eg_2}$ and  $M_{fg_2}$ have no finite type I direct summands with infinite dimensional center, the conclusion follows again from Proposition \ref {P2.2} and Lemma \ref {L2.3}.

Since $M_{eg_3}$  is a finite type I algebra then  $eg_3$ is a finite projection and  $M_{g_3}$ is a type I algebra. $M_{g_3}$ must have a finite dimensional center, because otherwise the center of $M_{fg_3}$ would be infinite dimensional as well, and since $fg_3$ is finite this would contradict the assumption that $M_f$ has no finite type I direct summands with infinite dimensional center. Since $(e+f)g_3$ is finite, we thus see that  $M_{(e+f)g_3}$ is a finite sum of matrix algebras. Therefore $ag_3\in M_{(e+f)g_3}$ is a positive combination of projections by the spectral theorem.

Finally, consider $ag_4$. Choose an integer $k> \frac{\|b\|}{\nu}$ and by   \cite [Lemma 3.2] {Hh77},
 decompose $fg_4$ into a sum $fg_4=\sum_{j=1}^kf_j$ of $k$ mutually orthogonal projections $f_j\sim fg_4$ that \textit{commute} with $dg_4$ and hence with $d$.
Then 
 $$
 ag_4= bg_4+dg_4= \sum_{j=1}^k\big(\frac{1}{k}beg_4+ df_j\big). 
 $$
 For every $1\le j\le  k$, $eg_4\perp f_j$, $eg_4 \prec  fg_4\sim f_j$, $f_j$ is properly infinite, $M_{f_j}$ has no  finite type I direct summands with infinite dimensional center, and  $$df_j=f_jd\ge \nu f_j\quad\text{with}\quad \nu > \| \frac{1}{k}beg_4\|.$$
Thus  by the first part of the proof   it follows that $\frac{1}{k}beg_4+ df_j$ is a  positive combinations of projections in $M$ for every $1\le j\le  k$   
  and hence so is $ag_4$. This concludes the proof.
 \ep

\bR{R:2.8}
\item [(i)] 
The first part of the proof  with the additional  assumption that  $\nu > \|b\|$  holds without changes also for  any C$^*$-algebra for which all positive locally invertible element  are positive combinations of projections. These algebras include, among others, all properly infinite simple $\sigma$-unital  C$^*$-algebras and their multiplier algebras. We will  focus on such C$^*$-algebras  in \cite {KNZC*10}. 
\item [(ii)] A key tool in this proof and also in the proof of  Theorem \ref {T4.3} below, is the fact that in a properly
infinite von Neumann algebra, given a selfadjoint operator (or, equivalently, given a masa) the identity  can be decomposed into the sum of two
equivalent projections commuting with the operator (equivalently, belonging to the masa). This result   was established by Halpern
in \cite [Lemma 3.2] {Hh77} and then obtained by different methods in the case of $\sigma$-finite algebras by Kadison \cite{Kr84}
(see also an extension by  Kaftal in \cite {Kv91}).
 \eR
Notice that the operator $a$ in Lemma \ref {L:2.7} satisfies the condition
\be{e:2.2} \exists~ \delta > 0 ~\text{ such that } ~\chi_a(0, \delta)\prec \chi_a[\delta,\infty).
\ee
It is easy to see, and it can also be obtained as a simple consequence of Lemma \ref {L:2.9} below, that if $M_{R_a}$ is
properly infinite, then condition (\ref {e:2.2}) is equivalent to
\be{e:2.3}
 \exists~ \delta > 0 ~\text{ such that } ~ \chi_a[\delta,\infty)\sim R_a.
\ee
In Theorem \ref {T:2.11}, we will see that condition (\ref{e:2.2}) is sufficient  for $a$ to be a finite sum of projections
 and is also necessary when $R_a$ is properly infinite.

First, we need a result  for which we could not find an explicit reference. It the generalization of the well know fact for $B(H)$ that if the supremum of a finite number of projections is infinite then at least one of the projections must be infinite.

\bL{L:2.9}
Assume that $M$ is properly infinite and that  $\bigvee_{j=1}^n p_j\sim I$ for some projections $p_j\in M$. Then there is a family of mutually orthogonal central projections $g_j$, some of which may be zero, with  $\sum_{j=1}^n g_j=I$ and for which $p_jg_j\sim g_j$ for every $j$.
\eL
\bp The proof is by induction on $n$. Let us  first prove the claim for $n=2$.
Assume first that $p_1\vee p_2=I$. Since $I$ is properly infinite, we can find $q\sim q^\perp\sim I$. By the comparison property of projections, there is a central projection $g$ such that \begin{align} (q\wedge p_1)g \prec (q^\perp\wedge p_1^\perp)g\label{e:2.4.0}\\
(q^\perp\wedge p_1^\perp)g^\perp\prec (q\wedge p_1)g^\perp.\label{e:2.5.0}\end{align}
Then
\begin{alignat*}{2}
g&\sim qg &(\text{since }q\sim I)\\
&=(q\wedge p_1)g+ (q- q\wedge p_1)g\\
&\prec  (q^\perp\wedge p_1^\perp)g +(q- q\wedge p_1)g &(\text{by (\ref {e:2.4.0}), since }q^\perp\wedge p_1^\perp ~\perp  q- q\wedge p_1)\\
&\sim (q^\perp\wedge p_1^\perp)g + (q\vee p_1-p_1)g \qquad &(\text{by Kaplanski's Parallelogram Law, since }q^\perp\wedge p_1^\perp ~\perp ~q\vee p_1-p_1 )\\
&= p_1^\perp g\\
&=(p_1\vee p_2-p_1)g\\
& \sim (p_2-  p_1\wedge p_2)g &(\text{by Kaplanski's Parallelogram Law})\\
& \le p_2g\\
&\le g.
\end{alignat*}
Thus $p_2g\sim g$. Similarly,  $p_1g^\perp\sim g^\perp$.

\

Next, consider the case when $p_1\vee p_2\sim I$, i.e.,  there is an isometry $w$ such that $w^*(p_1\vee p_2)w=I.$ Then $w^*p_1w \vee w^*p_2 w= I$, hence by the first part of the proof, there is a central projection $g$ such that $w^*p_1wg\sim g$ and $w^*p_2wg^\perp \sim g^\perp $. But then $p_1g\sim g$ and  $p_2g^\perp\sim g^\perp$ which concludes the case $n=2$.

\

Now assume that the property holds for $n-1$ and that
$
\bigvee_{j=1}^n p_j\sim I.
$ Then $p_1 \vee \Big(\bigvee_{j=2}^n p_j\Big)\sim I$ and hence by the result for $n=2$, there is a central projection $g_1$ for which
$p_1g_1\sim e_1$ and 
\be{e:2.6}
\Big(\bigvee_{j=2}^n p_j\Big) e_1^\perp \sim e_1^\perp.
\ee
If $e_1^\perp=0$, i.e., $e_1=I$, we choose $e_j=0$ for $j\ge 2$ and we are done. If $e_1^\perp\ne 0$, then by (\ref{e:2.6}),  $\bigvee_{j=2}^n \big(p_je_1^\perp\big) \sim e_1^\perp,$ hence by the induction hypothesis applied to the projection $\{p_je_1^\perp\}_2^n$ in the properly infinite algebra $M_{e_1^\perp}$, we obtain
a decomposition of $e_1^\perp= \sum _{j=2}^ne_j$ into mutually orthogonal central projections in $M_{e_1^\perp}$ for which $$p_je_j=p_je_1^\perp e_j\sim e_j.$$ Of course, the central projections 
$\{e_j\}_2^n$ thus found in the center of $M_{e_1^\perp}$ can be identified with central projections in $M$  (e.g., see \cite[Ch I Sect 2 Corollaire, Proposition 2]{Dj72}), which concludes the proof.
\ep

The following lemma is based on  and partially overlaps with the proof of \cite[Lemma 2.6]{Kv77}.

 \bL{L:2.10}
If  $a\in M^+$ and assume that $a= \sum_j \lambda_jp_j$ for a (finite or infinite) collection  of scalars  $\lambda_j>0$ and projections $p_j\in M$,  where the sum converges in the strong topology in case the collection is infinite. Assume furthermore that $\delta := \inf \lambda_j >0$.
Then for every $j$, $\chi_a(0, \delta) \prec R_a- p_j$ and $p_j\prec \chi_a[\delta, \infty).$ 
\eL
\bp If $a\ge \lambda p$ for some $\lambda >0$ and some projection $p\in M$, then 
\ba \|a^{1/2}\xi \|&\ge \lambda ^{1/2}\|\xi\| \hspace{-3cm}&&\forall \xi\in pH\\          
\phantom{ssssssssssssssssssss}\|a^{1/2}\xi\|&<  \lambda ^{1/2}\|\xi\| &&\forall 0\ne \xi\in \chi_{a^{1/2}}(0, \lambda^{1/2})H.\phantom
{sssssssssssssssssssssssssss}
\end{align*}
Since $ \chi_a(0, \lambda) = \chi_{a^{1/2}}(0, \lambda^{1/2})$, we see that  $ \chi_a(0, \lambda) H\cap pH=\{0\}$, i.e., $ \chi_a(0, \lambda)\wedge p =0$. But then, by  Kaplanski's Parallelogram Law,  $$  \chi_a(0, \lambda)= \chi_a(0, \lambda)-  \chi_a(0, \lambda)\wedge p\sim (R_a- p)-  (R_a- p)\wedge (R_a-\chi_a(0, \lambda))\le R_a- p $$ and thus
$$
 \chi_a(0, \lambda)\prec R_a- p.
$$
Also $$p= p - p\wedge  \chi_a(0, \lambda)\sim \big(R_a- \chi_a(0, \lambda)\big)-  \big(R_a- \chi_a(0, \lambda)\big)\wedge (R_a-p)\le R_a- \chi_a(0, \lambda)= \chi_a[\lambda, \infty),$$ and thus
$$p\prec \chi_a[\lambda, \infty).$$
Since $a \ge \lambda_jp_j$ for every $j$, we thus have
$$
\chi_a(0, \delta)\le  \chi_a(0, \lambda_j)\prec R_a-p_j \qquad\text{and}\qquad p_j\prec \chi_a[\lambda_j, \infty)\le \chi_a[\delta, \infty).$$
 \ep

 \bT{T:2.11}
 Let $M$ be a von Neumann algebra and let $a \in M^+$.
 \item[(i)] Assume that $M_{R_a}$ has no finite type I direct summands with infinite dimensional center and that    there is a $ \delta > 0$  such that  $\chi_a(0, \delta)\prec \chi_a[\delta,\infty)$. Then $a$ is a positive combination of projections in $M$.
 \item[(ii)] If $R_a$ is properly infinite, then $a$ is a positive combination of projections in $M$ if and only if  there is a $ \delta > 0$  such that  $\chi_a(0, \delta)\prec \chi_a[\delta,\infty)$. 
  \item[(iii)]  If  $M_{R_a}$  is a finite sum of finite factors or of  $\sigma$-finite type III factors, then $a$ is always a positive combination of projections in $M$.
 \eT
\bp
 \item[(i)]
Assume that $\chi_a(0,  \delta)\prec  \chi_a[\delta,\infty)$  for
some $\delta$ and set  \ba 
e&:=\chi_a(0, \delta) & f&:= \chi_a[\delta,\infty)\\
b&:= a\chi_a(0,  \delta) =ae &d& := a \chi_a[\delta,\infty)=af.
\end{align*}
Then 
$ a= b+c$ and $R_a=e+f.$
Since $c(e)\le c(f)$, it follows that $c(f)=c(R_a)$. Assume by contradiction that $M_f$ had a 
finite direct summand $M_{fg}$ of type I with infinite dimensional center for some central projection $g\le c(f)$. Then it would follow that $M_{g}$ is a type I algebra with infinite dimensional center and that $fg$ is a finite projection. But then $eg$ too is finite and hence so is $R_ag$. But then $M_{R_ag}$ too would be a finite type I algebra with infinite dimensional center against the assumption. Thus the hypotheses of Lemma \ref {L:2.7}  applied to $a=b+d$ are satisfied and hence $a$ is a positive combination of projections.
\item[(ii)] Since $R_a$ is properly infinite, then $M_{R_a}$ has no finite direct summands of any type, so $a$ is a positive combination of projections in $M$ by (i).  

Conversely, assume that $a$ is a positive combination of projections, i.e., $a= \sum_{j=1}^n \lambda_jp_j$ for some scalars $\lambda_j>0$ and projections $p_j$. To simplify notations, assume that $R_a=I$, i.e., $\bigvee_1^np_j=I$. Let $\delta= \min_j \lambda_j$.
By Lemma \ref {L:2.10}, $p_j\prec \chi_a[\delta, \infty)$ for all $1\le j\le n$ and 
by Lemma \ref {L:2.9}, there is a family of mutually orthogonal central projections $e_j$ with  $\sum_{j=1}^n e_j=I$
 for which $p_je_j\sim e_j$. Then $e_j\prec \chi_a[\delta, \infty)e_j$ for all $j$, and hence, $\chi_a[\delta, \infty)\sim I$. Thus (\ref{e:2.2}) holds.
\item[(iii)] Assume, without loss of generality, that $R_a=I$ and that $M$ is a factor. 

 If $M$ is of  type $I_n$, i.e., a matrix algebra, then the conclusion follows from the spectral theorem.

If $M$ is of type $II_1$ and $\tau$ is the canonical trace, then $\tau(\chi_a(0, \delta))\to 0$ for $\delta\to 0$. Choose a $\delta>0$ so that $\tau(\chi_a(0, \delta))\le \tau(\chi_a[\delta,  \infty))$. Then 
 $\chi_a(0, \delta))\prec \chi_a[\delta, \infty)$ and hence $a$ is  positive combination of projections by (i).

If  $M$ is a $\sigma$-finite type III factor,  then all nonzero
projections are equivalent,  and hence, the condition (\ref{e:2.2})
holds trivially. Again, the conclusion follows from (i). \ep

The following example shows that condition (\ref{e:2.2}) may not be necessary when $M$ is neither properly infinite nor it is a finite sum of finite factors.

\begin{example}\label{E:2.12}
Assume that $M= \sum_{j=1}^\infty M_j$ is an infinite direct sum of
finite factors $M_j$ such that  each $M_j$ contains three mutually
orthogonal and equivalent nonzero projections $e_j, f_j, g_j$.
Choose a sequence $(0,1) \ni \lambda_j \to 1$  and define $$a_j:=
(1-\lambda_j )e_j+ (1-\lambda_j )f_j+(1+2\lambda_j )g_j$$ and
$a:=\bigoplus _{j=1}^\infty a_j$. By Theorem \ref{T1.1}), each
$a_j$ is the sum of three projections,  and hence, so is $a$. In
particular, $a$ is a positive combination of projections. However,
for every $0<\delta< 1$, we have
\ba\chi_a(0, \delta)&= \bigoplus _{j=1}^\infty \chi_{a_j}(0, \delta)= \bigoplus _{1-\lambda_j< \delta} (e_j+f_j)\\
\chi_a[ \delta, \infty) &= \bigoplus _{j=1}^\infty \chi_{a_j}[ \delta, \infty) = \bigoplus _{j=1}^\infty g_j+
\bigoplus _{1-\lambda_j\ge \delta} (e_j+f_j).
\end{align*}
Since  $e_j+f_j\not \prec g_j$, it follows that $\chi_a(0, \delta)\not \prec \chi_a[ \delta, \infty).$
\end{example}
\section{Central essential spectra} \label{S:central}

Fong and Murphy  in \cite[Theorem 11] {FoMj85}  proved that if $a\in
B(H)^+$ ($H$  infinite dimensional and separable) and $R_a$ is infinite, then $a$ is a
positive combination of projections if and only if $a\not\in K(H)$,
i.e., if and only if $\|a\|_e>0$, where $\|a\|_e$ denotes the
essential norm. A natural extension of this result is to the case of
 $\sigma-$finite properly infinite von Neumann factors, by replacing  $K(H)$ by the
ideal of relative compact operators $J(M)$ generated by the finite projections of $M$.  This can be deduced
easily from Theorem \ref {T:2.11} and is a special case of Theorem
\ref {T3.4} below.

For von Neumann algebras that are not $\sigma-$finite, however, we need more general \textit{central ideals} and if the algebras are not factors, we need the corresponding center-valued version of the essential spectrum, the \textit{central essential spectrum}. The need for the latter is illustrated by the following example.

\begin{example}\label{E3.1}
Let $M=\bigoplus_{k=1}^\infty B(H)$, where $H$ is an infinite dimensional separable Hilbert space     and let  $a=\bigoplus_{k=1}^\infty a_k\in M^+$.  Assume that $R_a$ is properly infinite. Then $a$ is a positive combination of projections if and only if  for some  $\nu >0$
$$C_u(a):=\bigoplus_{k=1}^\infty \big(\|a_k\|_eI_k\big)~  \ge  ~\nu c(R_a).$$
Notice that $J(M)= \bigoplus_{k=1}^\infty K(H)$, hence $\|a\|_e= \sup_k \|a_k\|_e$. Thus in particular the condition that $\|a\|_e>0$ is not sufficient for $a$ to be a positive combination of projections.
\end{example}
\bp
For the necessity of the condition, let $a=\sum_{j=1}^n \lambda_jp_j$  for some $\lambda_j>0$ and projections $p_j= \bigoplus_{k=1}^\infty p_{jk}$.  Then
$a_k= \sum_{j=1}^n \lambda_jp_{jk}$ for every $k$. Let $\pi: B(H)\to B(H)/K(H)$ denote  the canonical quotient map. Then since
$\pi(a_k)= \sum_{j=1}^n \lambda_j\pi(p_{jk})$, it follows that
$\|a_k\|_e:= \|\pi(a_k)\|\ge \lambda_j\|\pi(p_{jk})\|$ for every $k$ and $j$. Now $ \|\pi(p_{jk})\|\in \{0,1\}$ and by the assumption that $R_a$ is properly infinite, all the summands $a_k$ are either $0$ or have infinite rank, i.e.,
$(R_a)_k= R_{a_k}= \bigvee_{j=1}^n p_{jk}$ is either $0$ or an infinite projection. Thus  if $a_k\ne 0$, then there is at least one index $j$ for which $p_{jk}$ is infinite, in which case $\pi(p_{jk})\ne 0$ and hence
$\|\pi(p_{jk})\|=1$. But then,  $\|a_k\|_e\ge  \min \lambda_j$. Let $\nu: =\min \lambda_j.$ Then $\nu > 0$ and $\bigoplus_{k=1}^\infty \big(\|a_k\|_eI_k\big) \ge \nu c(R_a).$

For the sufficiency of the condition, notice that for any $0< \delta < \nu$  and for every $k$ for which $a_k\ne 0$ and hence $\|a_k\|_e\ge \nu$ it follows that $\chi_{a_k}[\delta, \infty)\sim I_k$ and hence
$$
\chi_a[\delta, \infty)=\bigoplus_{k}\chi_{a_k}[\delta, \infty)\sim \bigoplus_{a_k\ne 0}I_k\sim R_a.
$$
Thus the condition (\ref{e:2.2}) is satisfied and hence by Theorem \ref {T:2.11}, $a$ is a positive combination of projections.  

\ep

The  operator  $C_u(a)=\bigoplus_{k=1}^\infty\big( \|a_k\|_eI_k\big)$  is largest element in what is called the \textit{central essential spectrum} of $a$, which in the case of a direct sum of factors is easy to describe as:
$$Z\text{-Sp}(a):=\big \{\bigoplus_{k=1}^\infty \big(\mu_kI_k\big)\mid \mu_k\in \sigma_e(a_k)\big\}.$$
For arbitrary von Neumann algebras, central essential spectra were introduced   in the study of central ideals by Halpern \cite{Hh72}. Example \ref {E3.1} suggest a reformulation of  the condition (\ref {e:2.2}) in  Theorem \ref{T:2.11} (ii) in terms of central ideals and central essential spectra.

For the reader's convenience we summarize  the notions from \cite {Hh72} concerning  central ideals and central essential spectra
that   will be need  here. The reader interested only in $\sigma$-finite factors can skip the remainder of  this section.

An ideal $\mathscr J$ in a von Neumann algebra $M$  is called  \textit{central}  if whenever $\{g_\gamma\}_{\gamma\in \Gamma}$ is a family of mutually orthogonal central projections $g_\gamma$  and $\{a_\gamma\}_{\gamma\in \Gamma}$ is a norm-bounded family of elements of $\mathscr J$, then $a:=\sum_{\gamma\in \Gamma}a_\gamma g_\gamma$ belongs to $\mathscr J$.

Denote by $Z:=M\cap M'$ the center of  the von Neumann algebra $M$. For every properly infinite projection $p\in M$ and every central projection $f \ge c(p)$, (recall that $c(p)$ denotes the central support of $p$,) let $\mathscr J_f(p)$ be the norm closed two sided ideal of $M$ whose set of projections is
$$\{q\in M \mid  q\le f \text{ and }pg \prec  qg \text{ for some projection } g\in Z \Rightarrow pg=0.\}$$
Thus by definition $$\mathscr J_f(p)= \mathscr J_{c(p)}(p)\oplus M_{f-c(p)}$$ and a projection $q$ belongs to $\mathscr J_{c(p)}(p)$ if and only if $q\prec p$ and $qg\not \sim pg$ for all central projections $0\ne g\le c(p).$
Then $\mathscr J_f(p)$  is central and every  central ideal $\mathscr J$ has this form.
If we further assume that $\overline{\mathscr J_f(p)}^w=Mf$, then $f$ is unique and $p$ is unique up to Murray-von Neumann equivalence (\cite[Theorem 2.4, Propostion 2.5] {Hh72}).

To simplify notations, identify  $Z$ with $C(X)$, with  $X$ a compact Hausdorff extremely disconnected space. $X$ is called the spectrum of the abelian algebra $Z$ and  is identified by the Gelfand theorem with the collection of the closed two-sided maximal ideals of $Z$. For every $\zeta\in X$ let $[\zeta]$ denote the closed two-sided ideal of $M$ generated by $\zeta$. For every closed two-sided ideal $\mathscr J$ of $M$,  $\mathscr J+[\zeta]$ is also a closed two-sided ideal of $M$ and $\pi_{\mathscr J+[\zeta]}$ denotes the canonical quotient map $M\to M/(\mathscr J+[\zeta])$.

Given a central ideal $\mathscr J$ and an element $a\in M$, the central essential spectrum of $a$ is defined as the collection
$$
Z\text{-Sp}(a):= \{z\in Z \mid z(\zeta)\in \text{Sp}\big( \pi_{\mathscr J+[\zeta]}(a)\big) \text{ for all  } \zeta\in X\}.
$$
Then   $Z\text{-Sp}(a)$ is always a non-empty subset of the center (more precisely, of the center of $M_{c(R_a)}$) and  is closed in the strong operator topology \cite [Theorem 3.5, Proposition 3.9] {Hh72}.
If  $a=a^*$, then $Z\text{-Sp}(a)$  is selfadjoint and has \cite [Proposition 3.12]{Hh72} a maximal element  $C_u(a)\in Z\text{-Sp}(a)$
given  for all
$\zeta \in X$  by
\ba
C_u(a)(\zeta):&= \max \text{Sp}\big( \pi_{\mathscr J+[\zeta]} (a)\big).
\end{align*}
When $a$ is positive, $C_u(a)$ takes the role of ``central essential norm".
Since we consider only one central ideal at a time, to simplify notations   we will not mark explicitly the dependence on $\mathscr J$ of $Z\text{-Sp}(a)$  or of the maximal
element  $C_u(a).$

Recall the well-know property in $B(H)$ ($H$ infinite dimensional  separable): if $a\in B(H)^+$, then $\|a\|_e\ge \nu$   if and only if for every $\epsilon >0$ the spectral projection  $\chi_a[\nu-\epsilon, \infty)$ is infinite (equivalently, $\chi_{a}[\nu-\epsilon, \infty)\sim I.$) It  is implicit in \cite {Hh72} that  this characterization extends to the central essential spectrum  in the context of properly infinite von Neumann algebras, but for the readers' convenience we make the connection explicit, listing only the facts that we need for  this article.

From now on,  fix  an $a\in M^+$ with properly infinite range projection $R_a$ and  set $\mathscr J:= \mathscr J_{c(R_a)}(R_a)$.
Notice that the  operator $\bigoplus_{k=1}^\infty \big(\|a_k\|_eI_k\big)$ in Example \ref {E3.1}  is indeed  the element $C_u(a)$ relative to this ideal $\mathscr J$.

\bL{L3.2}
If $z\in Z\text{-Sp}(a)$, $\epsilon > 0$ and $p:= \chi_{a-z}[-\epsilon, \epsilon]$. Then $R_a\prec p$.
\eL
\bp
To simplify notations, assume that $c(R_a)=I$.  Reasoning as in \cite [Proposition 2.9]{Hh72}, let  $f$ be the maximal  central projection  for which $R_af\prec pf$.   If $R_ag\prec pf^\perp g$ for some central projection $g$, then also $R_af^\perp g\prec pf^\perp g$. By the  maximality of $f$, it follows that $f^\perp g=0$, and hence that $R_ag=0$.  By definition,  $pf^\perp\in \mathscr J$. By \cite [Proposition 3.13]{Hh72}, $ f^\perp=0$ and thus $R_a\prec p$.
\ep
\bP{P3.3}
Let  $\nu > 0$. Then $C_u(a)\ge \nu c(R_a)$ if and only if  $\chi_a[\nu-\epsilon, \infty)\sim R_a$ for all $0<\epsilon < \nu.$
\eP
\bp
Without loss of generality, assume that $c(R_a)=I$.  For a fixed  $0<\epsilon < \nu$, to simplify notations  let $z:= C_u(a)$ and $p: =  \chi_{a-z}[-\epsilon,~\epsilon]$.
Assume first that $ z= C_u(a)\ge \nu c(R_a)$.  Since $p$ commutes with $a-z$, it commutes also with $a$.  Then
$$(\nu-\epsilon) p \le zp-\epsilon p \le ap\le zp+\epsilon p\le (\|z\|+\epsilon)p.
$$
Hence  $$p\le \chi_a[ \nu-\epsilon,~ \|z\|+\epsilon]\le \chi_a[\nu-\epsilon, \infty)\le \chi_a(0, \infty)=R_a.$$
But $R_a\prec p$ by Lemma \ref {L3.2}, hence  $R_a\sim \chi_a[\nu-\epsilon, \infty).$

Assume now that for every $0<\epsilon < \nu$ we have  $r:=\chi_a[\nu-\epsilon, \infty)\sim R_a$. Then $a\ge ar\ge (\nu-\epsilon)r. $ Thus for every $\zeta \in X$ we have
$$
\pi_{J+[\zeta]}(a) \ge (\nu-\epsilon)\pi_{\mathscr J+[\zeta]}(r),
$$
and hence
$$C_u(a)(\zeta)=\max\{\text{Sp}\big( \pi_{\mathscr J+[\zeta]}(a)\big)\}= \|\pi_{\mathscr J+[\zeta]}(a)\|\ge (\nu-\epsilon)\|\pi_{\mathscr J+[\zeta]}(r)\|.$$
Since $\pi_{\mathscr J+[\zeta]}(r)$ is a projection, $\|\pi_{\mathscr J+[\zeta]}(r)\|$ is either $0$ or $1$.
Reasoning by contradiction, assume that for some $\zeta_o \in X$, we have $\|\pi_{\mathscr J+[\zeta_o]}(r)\|= 0$ .
The function $$X\ni \zeta\to \|\pi_{\mathscr J+[\zeta]}(r)\|\in \mathbb R$$
is continuous by \cite [Theorem 3.2]{Hh72}, and hence $\|\pi_{\mathscr J+[\zeta]}(r)\|$ vanishes on a closed set containing $\zeta_o$. This set is clopen because the spectrum of any abelian von Neumann algebra is extremely disconnected,  thus the characteristic function of this set is also continuous. Hence that characteristic function is identified with a  projection $g$ in the center of $M$. Then by \cite [Lemma 3.1] {Hh72}, $rg\in \mathscr J.$ Since $r\sim R_a$,  reasoning as in \cite [Proposition 2.9]{Hh72}, it follows that $g=0$, a contradiction. Thus
$\|\pi_{\mathscr J+[\zeta]}(r)\|= 1$ for all $\zeta$, i.e.,  $C_u(a)(\zeta)\ge \nu-\epsilon$ for all $\zeta$. Thus $C_u(a)\ge (\nu-\epsilon)I$ and hence $C_u(a)\ge \nu I$.
\ep

Notice that the existence of $\nu>0$ such that $\chi_a[\nu-\epsilon, \infty)\sim R_a$ for all $0<\epsilon < \nu$ is obviously equivalent to  the existence of $\delta >0$ such that $\chi_a[\delta, \infty)\sim R_a$
(condition (\ref {e:2.3})), which in turns as we have already remarked, is  equivalent to condition (\ref{e:2.2}). In other words, we can reformulate condition (ii) in Theorem  \ref {T:2.11} in terms of central essential spectra as follows.

\bT{T3.4}
Let $a\in M^+$ have  properly infinite range projection $R_a$ and let $J= J_{c(R_a)}( R_a)$. Then the following are equivalent:
\item[(i)] $a$ is a positive combination of projections in $M$.
\item[(ii)] There is a $\delta>0$  for which $\chi_a(0,\delta)\prec \chi_a[\delta,
\infty)$.
\item [(iii)] There is a $\delta>0$  for which $\chi_a[\delta, \infty)\sim
R_a$.
\item[(iv)] $C_u(a)$ is locally invertible, i.e., there is a $\nu > 0$ for which $C_u(a)\ge \nu c(R_a)$.
\eT

In the case when $M$  is a properly infinite $\sigma-$finite factor, and $R_a$ is infinite, then $C_u(a)=\|a\|_eI$ where $\|\cdot\|_e$ is the essential norm relative to the relative compact ideal $J(M)$. Thus condition (iv) is equivalent to asking that $a\not \in J(M)$. This the generalization of the Fong and Murphy characterization in $B(H)$ in \cite[Theorem 11] {FoMj85}.

\bC{C3.5}
Let $M$ be a properly infinite $\sigma-$finite factor and let $a\in M^+$ have a  properly infinite range projection $R_a$. Then $a$ is positive combination of projections if and only if $a\not\in J(M)$.
\eC

\section{Finite sums of projections in properly infinite algebras}\label{ess norm >1}
In this section we prove, for general properly infinite von Neumann
algebras  an analog of the $B(H)$ result in \cite {WpCm09} that
positive operator with essential spectrum larger than 1 are finite
sums of projections. Our approach, which is different from the one by Choi and Wu,
is to  decompose the operator into an infinite sum of blocks, where each block is a sum of not more than a fixed number of projections and where the
blocks separated by more than one index are orthogonal. The
following simple folklore lemma shows how to then reassemble this
infinite sum into a finite sum of projections.

\bL{L4.1}
Let $\{b_j\}$ be a sequence of positive elements of $M$ such that
\item[(i)] $R_{b_i}R_{b_j}=0$ for all $|i-j|>1$,
\item[(ii)] $b_j$ is the sum of $N_j$ projections in $M,$  and
\item[(iii)] $N:=\sup N_j< \infty$.\\
Then $b:= \sum_{j=1}^\infty b_j$ is the sum of $2N$ projections in $M.$
\eL
\bp
Let $$b_o:=\sum_{j=1}^\infty b_{2j-1}\qquad\text{and}\qquad b_e:=\sum_{j=1}^\infty b_{2j}.$$ By hypothesis, the range projections of the elements $b_{2j-1}$ are mutually orthogonal, and so are those of the elements $b_{2j}$. Since $\sup_j \|b_j\|\le N$,  both series converge in the strong operator topology  and hence $b=b_o+b_e.$  Set $b_j=\sum_{k=1}^N q_{j,k}$ where $q_{j,k}$ are projections in $M$ (some of these projections can vanish.)
For each $1\le k\le N$, let 
$$q_{o,k}:=\sum_{j=1}^\infty q_{2j-1,k}\qquad\text{and}\qquad 
q_{e,k}:=\sum_{j=1}^\infty q_{2j,k}.
$$
Both series being the sum of mutually orthogonal projections converge to a projection. Thus
$b_o= \sum_{k=1}^Nq_{o,k}$ and $b_e= \sum_{k=1}^Nq_{e,k}$ are each the sum of $N$ projections in $M$ and hence $b=b_o+b_e$  is the sum of $2N$ projections in $M.$
\ep

The following is one of the key lemmas in this paper. The core of the argument  is based on Fillmore's characterization  of finite sums of finite projections (see Theorem \ref{T1.1}).

\bL{L4.2}
 Let $M$ be a von Neumann algebra and $e$ and  $f$ be orthogonal projections in $M$ such that  $e\prec f$ and $f$ is properly infinite. Let $a=b e +\alpha f$ where $ \alpha>1$ and $ b=be=eb \ge 0$,  Then $a$ is the sum of  finitely many  projections in $M$.
\eL
\bp
We first reduce the problem to the case when $b$ is a scalar multiple of $e$ and $f=\sum _{k=1}^\infty e_k$ is an infinite sum of projections equivalent to $e$.

Decompose $f=f'+f''$ into a direct sum of two equivalent
projections, $f'\sim f''\sim f$. Since $f$ and hence $f'$ are properly infinite, $M_{f'}$ has no direct finite summands. Thus Lemma 
\ref {L:2.7} applies and $b+ \alpha f'$ is a positive combination of
projections. Say $$b+ \alpha f' =\sum_{j=1}^n \beta_jp_j$$ for some $n\in \mathbb N$, $\beta_j>0$, and projections $p_j\in M.$
Further decompose the properly infinite projection $f''= \sum_{j=1}^np_j$ into a sum of $n$ mutually orthogonal 
equivalent projections $f_j\sim f''\sim f$. Then
$$a= \sum_{j=1}^n \big(\beta_jp_j +\alpha f_j\big),$$
where $p_j\le e+f' \perp f''$ and hence $p_j\perp f_j$ for all $j$. Furthermore, 
$p_j\prec f\sim f_j $ and $f_j$ is properly infinite. Thus after this first part of the reduction, we can assume, without loss of
generality, that $a = \beta e +\alpha f$ for some $\beta >0$.

Next,  since the projections $fc(e)$ and $f(c(f)-c(e))$ either vanish or are properly infinite, we can further decompose them  into two infinite sums of mutually orthogonal  projections
$e'_k\sim e$ and $e''_k\sim f(c(f)-c(e)) $, that is, 

$$fc(e)= \sum_{k=1}^\infty e'_k\quad\text{and} \quad f(c(f)-c(e))= \sum_{k=1}^\infty e''_k.$$ Then $$
a=\beta e+ \alpha f(c(e))+ \alpha f(c(f)-c(e))
=\beta e+  \alpha \sum_{k=1}^\infty e'_k + \alpha \sum_{k=1}^\infty e''_k.$$
Both summand have the form $\beta e_0+ \alpha \sum_{k=1}^\infty e_k$.
where the projections $\{e_k\}_0^\infty$ are mutually orthogonal and equivalent, $\beta \ge 0$, and $\alpha >1.$
Thus this completes the reduction. Explicitly, we assume without loss of generality  that
$$ a=\beta e_0+ \alpha \sum_{k=1}^\infty e_k,\quad e_k\sim e_o, \quad\beta \ge 0, ~\alpha >1.
$$

The equivalence of the projections $\{e_k\}_0^\infty$, defines an embedding of $B(H)$ in $M$ ($H$ separable infinite dimensional)  under which the projections $\{e_k\}_0^\infty$ correspond to rank-one projections in $B(H)$. Thus to simplify notations,
assume henceforth that $M=B(H)$ and that  all the projections $e_k$ have rank one.

Since
$$a= (\beta -\lfloor\beta\rfloor) e_o+ (\alpha- \lceil\alpha-2\rceil) \sum_{k=1}^\infty e_k
 +\lfloor\beta\rfloor e_o+  \lceil\alpha-2\rceil \sum_{k=1}^\infty e_k
 $$  and
 $\lfloor\beta\rfloor e_o+  \lceil\alpha-2\rceil \sum_{k=1}^\infty e_k$ are already
 finite sums of projections, we can assume, without loss of generality, that $1<\alpha\le 2$ and that $0\le \beta <1$.
Set $\beta_0:=\beta$,  $n_0=0$ and $N:=  \lfloor\frac{\alpha^2}{\alpha-1}\rfloor.$

Let
\begin{align*}
n_1&:= \Big\lceil \frac{1-\beta_0}{\alpha-1}\Big\rceil \\
\beta_1&:= \beta_0+ n_1\alpha- \lfloor   \beta_0+ n_1\alpha \rfloor\\
a_1&:=\beta_0e_{n_0}+ \alpha\sum_{j=n_0+1}^{n_1-1}e_j+ (\alpha-\beta_1)e_{n_1}.
\end{align*}
Here and in the sequel, we adopt the convention of dropping  any sum where the upper index of summation strictly less than the lower index of summation. So,  if $n_1=1$, (which occurs when $1-\beta\le \alpha -1$),  then the above formula just reads $a_1= \beta_0e_{0}+ (\alpha-\beta_1)e_{1}$.

By definition,  $0\le \beta_1<1<\alpha$ and thus $a_1$ is a positive finite rank operator with
\ba
\tr(a_1)&=\beta_0 + (n_1-1)\alpha+ \alpha- \beta_1= \lfloor   \beta_0+ n_1\alpha \rfloor\in \mathbb N.
\\
 R_{a_1}&=\begin{cases} e_{n_0}+  \sum_{j=1}^{n_1}e_j\quad &\text{if } \beta_0\ne 0\\
\sum_{j=1}^{n_1}e_j\quad &\text{if } \beta_0= 0\end{cases}\\
 \intertext{and hence}
 \text{rank}(a_1)&=\begin{cases} n_1+1\hspace{1.55cm} &\text{if } \beta_0\ne 0\\
 n_1&\text{if } \beta_0= 0\end{cases}\hspace{1cm} \le \ \ n_1+1\\
\end{align*}

By definition of $ \lceil \cdot \rceil $,
$$
\frac{1-\beta_0}{\alpha-1}\le n_1< \frac{1-\beta_0}{\alpha-1}+1= \frac{\alpha-\beta_0}{\alpha-1},
$$
hence
$$ n_1+1 \le \beta_0+ n_1\alpha <  \beta_0+\frac{(\alpha-\beta_0)\alpha}{\alpha-1} = \frac{\alpha^2-\beta_0}{\alpha-1}\le \frac{\alpha^2}{\alpha-1}.
$$
From the first inequality,   it follows that 
$$\text{rank}(a_1)\le n_1+1 \le \lfloor  \beta_0+ n_1\alpha \rfloor = \tr(a_1). $$ From the remaining chain of inequalities it follows that
$$\tr(a_1)= \lfloor  \beta_0+ n_1\alpha \rfloor  \le \lfloor \frac{\alpha^2}{\alpha-1} \rfloor  = N.$$ Since 
$
\text{rank}(a_1) \le  \tr(a_1)\in \mathbb N,
$
 Theorem \ref{T1.1} applies and  we conclude that $a_1$ is the sum of $\tr(a_1)\le N$ rank-one projections.

 Now
$$a-a_1= \beta_0e_{n_0}+ \alpha\sum_{j=n_0+1}^{\infty}e_j-\Big(\beta_0e_{n_0}+ \alpha\sum_{j=n_0+1}^{n_1-1}e_j+ (\alpha-\beta_1)e_{n_1}\Big)= \beta_1 e_{n_1} + \alpha\sum_{j=n_1+1}^{\infty}e_j
$$
has the same form as the beginning operator $a$ and with the same constant $\alpha$. Thus repeating the construction, let
\begin{align*}
n_2&:=n_1+  \Big\lceil \frac{1-\beta_1}{\alpha-1}\Big\rceil  \\
\beta_2&:= \beta_1+ (n_2-n_1)\alpha- \lfloor \beta_1+ (n_2-n_1)\alpha\rfloor\\
a_2&:=\beta_1e_{n_1}+ \alpha\sum_{j=n_1+1}^{n_2-1}e_j+ (\alpha-\beta_2)e_{n_2}.
\end{align*}
Then by the same argument as above, $a_2$ is a positive operator with
\ba \text{rank }(a_2)&=\begin{cases}n_2-n_1+1\quad &\text{if } \beta_1\ne 0 \\
n_2-n_1\quad &\text{if } \beta_1= 0
\end{cases} \  \ \le \  \ n_2-n_1+1\\
 \tr(a_2)&=  \lfloor \beta_1+ (n_2-n_1)\alpha\rfloor \  \in  \ \mathbb N
 \end{align*}
 and $$\text{rank}(a_2)\le \tr(a_2)\le N.$$
So, again, by Theorem \ref{T1.1}), $a_2$ is a sum of not more than $N$ projections.
Notice that $$R_{a_1}R_{a_2}= \begin{cases} 0&\beta_1=0\\
e_{n_1}&\beta_1\ne 0,\end{cases}$$
 thus in general $R_{a_1}R_{a_2}$ does not vanish. 
We can iterate this process to construct for each $k$ the positive operator
$$a_k=\beta_{k-1}e_{n_{k-1}}+ \alpha\sum_{j=n_{k-1}+1}^{n_k-1}e_j+ (\alpha-\beta_k)e_{n_k}$$
which is the sum of at most $N$ projections. By construction, $a= \sum_{j=1}^\infty a_j$.
An immediate but key observation  is that already $R_{a_1}R_{a_3}=0$ and in general  $R_{a_j}R_{a_{j'}}=0$ for $|j-j'|>1$.
Thus by Lemma \ref {L4.1},
$a$ is a sum of finitely many projections.
\ep

Now we proceed to one of our main results.

\bT{T4.3} Let $M$ be a von Neumann algebra. Assume that
  $a\in M^+$ has a properly infinite range projection $R_a$
and $C_u(a)\ge \nu c(R_a)$ for some $\nu>1$. Then $a$ is a finite
sum of projections in $M$. \eT

 \begin{proof}
Let $1< \alpha < \nu$. By Proposition \ref {P3.3}, $\chi_a[\alpha, \infty)\sim R_a$. By \cite [Lemma 3.2] {Hh77}, $\chi_a[\alpha, \infty)=p+q$ for some  projections $p$ and $q$ that commute with $a$ and for which $p\sim q\sim R_a.$ Then in particular, $$ap = a\chi_a[\alpha, \infty)p\ge\alpha \chi_a[\alpha, \infty)p= \alpha p$$ and similarly, $aq\ge \alpha q.$
 Let
\begin{align*}
 a_1&:=a(R_a-p-q)+  (a- \alpha I)q+  \alpha p\\
 a_2&:=(a- \alpha I)p+  \alpha q
  \end{align*}
so that $a= a_1+a_2.$
Since $p$ and $q$  and hence $R_a-p-q$ commute with $a$, it follows that $a_1$ and $a_2$ are both the direct sum of positive operators. Since $p\sim q\sim R_a$, it follows that $R_a-p\prec p$ and $p\prec q$, thus both $a_1$ and $a_2$ satisfy Lemma \ref {L4.2} and hence  are both finite sums of projections.  Thus $a$ too is a finite sum of projections.
 \ep
As already mentioned, in the $\sigma-$finite factor case, the central essential norm reduces to the  (scalar) essential norm relative to the ideal $J(M)$. 
\bC{C4.4} Assume that $M$ is an infinite $\sigma-$finite factor and $a\in M^+$. A sufficient condition for $a$ to be a finite sum of projections in $M$ is that
\item [(i)] $\|a\|_e>1$ when $M$ is of type I$_\infty$ (usual essential norm of $B(H)$;) 
\item [(ii)] $\|a\|_e>1$ when $M$ is of type II$_\infty$  (essential norm relative to the ideal  $ J(M)$;)
\item [(iii)]  $\|a\|>1$ when $M$ is of type III (operator norm).
\eC
The following examples show that we cannot relax the condition that $C_u(a)$ is bounded away from $I$.

First, notice that if $a\in B(H)^+$ and $R_a$ is infinite, then the condition $\|a\|_e\ge 1$ is necessary
 for $a$ to be a finite sum of projections, since then at least one of the projections must be infinite, and hence, has essential norm 1. However, the following simple example shows that  the condition $\|a\|_e\ge 1$ is not sufficient. 
 
\begin{example}\label{E4.5}
Let $a=I +k$ with $k\in K(H)^+$,  $0< \tr (k)< \infty$ but $\tr (k)\not \in \mathbb N$. Then  $\|a\|_e=1$. By \cite [Theorem 17]{KNZW*09}, $a$  is not an infinite sum of projections, hence cannot be a finite sum of projections. In fact even if $\tr (k) \in \mathbb N$, if $k$ has infinite rank, then $a$ is not a finite sum of projections by Corollary \ref {C5.9}  below.
\end{example}

 Next notice that even the condition  $C_u(a)(\zeta) > 1$ for all $\zeta\in X$ is not sufficient.

\begin{example}\label{E4.6}
Let $H$ be an infinite dimensional separable Hilbert space, let $M:=\bigoplus_{n=1}^\infty B(H)$ and let $a:= \bigoplus_{n=1}^\infty (1+\frac{1}{n})I$. Then $C_u(a)=a\ge I$ but it is not bounded away from $I$.  By \cite [Theorem 3]{KRS02}, for $n\ge 5$, $(1+\frac{1}{n})I$ is the sum of not less than $n+1$ projections, and thus, $a$ cannot be the sum of finitely many projections in $M$.
\end{example}

 \section{Necessary conditions}\label{nec cond}
First, we extend  to general von Neumann algebras a result that we have obtained  for $\sigma$-finite factors and strong sums of projections in \cite [Proposition 3.1]{KNZW*09}.  For the $B(H)$ case that result is equivalent to \cite [Theorem 1.2] {WpCm09}.  Denote by $g$ the maximal finite central projection of $M$ and by $\Phi$ be the canonical center-valued trace on $M_g$.

\bP{P5.1}
 Let $a\in M^+$ and $n\in\mathbb N$. Then the following conditions are equivalent:
\item [(i)]
There is a partial isometry $v\in M$ with $v^*v=R_a$ and a decomposition of the identity $I=\sum_{j=1}^nq_j$ into n mutually orthogonal nonzero projections $q_j\in M$  for which $q_jvav^*q_j=q_j$ for $1\le j\le n$.
\item [(ii)] $a=\sum_{j=1}^np_j$ is the sum of $n$ nonzero projections $p_j\in M$   and there is a decomposition of the identity $I=\sum_{j=1}^nq_j$ into n mutually orthogonal nonzero projections $q_j\in M$ for which $p_j\sim q_j$ for $1\le j\le n$.
\item [(iii)] $a=\sum_{j=1}^np_j$ is the sum of $n$ nonzero projections $p_j\in M$, $\Phi(ag)= g$, and $R_ag^\perp \sim g^\perp$.
\eP
\bp
(i) $ \Leftrightarrow$ (ii). Basically, the same proof as in \cite [Proposition 3.1]{KNZW*09}, but restricted to finite sums. For the readers' convenience we briefly reproduce it here. Assume that (i) holds. For every $1\le j \le n$, let $w_j:=q_jva^{1/2}$. Then  $w_jw_j^*=q_jva^{1/2}a^{1/2}v^*q_j =q_j$ and hence $p_j:= w_j^*w_j$ is a projection  equivalent to $q_j$. Then 
$$\sum_{j=1}^np_j=\sum_{j=1}^na^{1/2}v^*q_jva^{1/2} =a^{1/2}v^* \Big( \sum_{j=1}^n q_j \Big) va^{1/2}=  a^{1/2}v^*va^{1/2}= a^{1/2}R_aa^{1/2}=a.$$

Conversely, assume (ii) holds. Let $w_j\in M$ be partial isometries for which $w_jw_j^*=q_j$ and $w^*_jw_j=p_j$ and let $b:=\sum_{j=1}^n w_j$. Then $b^*b= \sum_{i,j=1}^nw_i^*w_j=\sum_{j=1}^np_j =a$. Let $b=va^{1/2}$ be the polar decomposition of $b$. then $v^*v=R_a$ and $q_jva^{1/2}=q_jb= \sum_{i=1}^nq_jw_i= w_j$, hence $q_jvav^*q_j=q_j$.

(ii)$ \Leftrightarrow$ (iii) It is enough to  prove the claim separately for the cases when $M$ is properly infinite  and when it is finite.

Assume first that $M$ is properly infinite (i.e., $g=0$) and assume that (ii) holds. By Lemma \ref {L:2.9}, there is a decomposition of the identity $I=\sum_{j=1}^ng_j$ into central projections $g_j$ such that $q_jg_j\sim g_j$. Then for every $j$
$$
R_ag_j= \Bigg(\bigvee_{i=1}^np_i\Bigg)g_j \ge p_jg_j\sim q_jg_j\sim e_j.
$$
Thus $R_ag_j\sim g_j$ for all $j$ and hence $R_a\sim I$.
Assume next that (iii) holds, i.e.,  $R_a\sim I$.
Again by Lemma \ref {L:2.9}, there is a decomposition of the identity $I=\sum_{j=1}^ng_j$ into central projections $g_j$ such that $p_jg_j\sim g_j$. By the assumption that $M$ is properly infinite, we can decompose each $g_j=\sum_{i=1}^ne_{ij}$ into a sum of mutually orthogonal and equivalent projections $e_{ij}\sim g_{j}$. Then for each $i \ne  j$, $p_ig_j\prec e_{ij}$, thus we can choose projections   $p_{ij}$ so that $p_ig_j\sim p_{ij}\le e_{ij}$. Set  $p_{jj}:=   e_{jj}+\sum_{i\ne j}(e_{ij}- p_{ij})$. Then
$$
p_jg_j\sim g_{j}\sim e_{jj}\le p_{jj} \le  g_j,
$$
hence also $p_jg_j\sim p_{jj}$. The projections $p_{ij}$ are mutually orthogonal, $g_j=\sum_{i=1}^np_{ij} $ for all  $j$, and hence $\sum_{i,j=1}^np_{ij} =I$. Let
$q_i:=\sum_{j=1}^np_{ij} $. Thus $\sum_{i=1}^nq_i=I$. Furthermore,
$q_ig_j= p_{ij}\sim p_ie_j$ for all $j$ and hence $q_i\sim p_i$ for all $i$, which proves (ii).

Assume now that $M$ is finite (i.e., $g=I$) and that (ii) holds. Then
$$\Phi(a) = \Phi\Big(\sum_{j=1}^np_j\Big)=\sum_{j=1}^n \Phi(p_j)= \sum_{j=1}^n \Phi(q_j)= \Phi(I)=I,$$
which proves (iii). 
Finally, assume that (iii) holds, i.e., $I=\Phi(a)= \sum _{j=1}^n\Phi(p_j)$. Set $q_1:=p_1$.
Then $\Phi(p_2)\le I- \Phi(p_1)=\Phi(I-p_1)$, hence $p_2\prec I-q_1$ and thus there is a projection $q_2$ with  $p_2\sim q_2\le I-p_1$. Proceeding inductively, we find mutually orthogonal projections $q_j\sim p_j$. Then 
$$
\Phi(I-\sum_{j=1}^nq_j) = I-  \sum _{j=1}^n\Phi(q_j)= I-  \sum _{j=1}^n\Phi(p_j) =0,$$
whence $\sum_{j=1}^nq_j=I$.
\ep

We need also the following result about two-sided not necessarily closed ideals.

\bL{L5.2}
Let $M$ be a von Neumann algebra let $I= \sum_{j=1}^nq_j$ for some finite collection of mutually orthogonal projections $q_j\in M$, let  $\Psi(x)$ be the map on $M$ defined by $ \Psi(x):= \sum_{j=1}^n q_jxq_j
\in M.$ Then 
\item [(i)] $\Psi$  is  a faithful, normal, and trace-preserving conditional expectation onto  $\bigoplus _{j=1}^nq_jMq_j$. In particular, it is linear and positive.
\item [(ii)] If $J$ is a  two-sided ideal of $M$ (not necessarily closed) and  $a\in M^+$, then $a\in J$ if and only if $\Psi(a)\in J$.
\eL
\bp
\item [(i)] Well known.
\item [(ii)] If $a\in J$, then obviously, $ \Psi(a)= \sum_{j=1}^n q_jaq_j
\in J.$  Assume that $\Psi(a)\in J$. Since $J$ is hereditary, it follows that  $q_jaq_j\in J$  for every $j$.

Let $$J^{1/2}:= \text{span} \{ x\in M^+\mid x^2\in J\}.$$
It is well known that  $J^{1/2}$ is also an ideal. (This fact can be easily verified using the characterization of the positive part of an ideal as a hereditary and additive collection of positive operators, see for instance \cite [3.21] {SsZl79}.) For every $1\le i,j\le n$,
$$|q_i a^{1/2}q_j|^2= \big(q_i a^{1/2}q_j\big)^*\big(q_i a^{1/2}q_j\big)= q_ja^{1/2}q_ia^{1/2}q_j\le q_jaq_j \in J,
$$
 whence $|q_i a^{1/2}q_j|\in J^{1/2}$, and hence, $q_i aq_j \in J^{1/2}$. As a consequence,
 $$
a^{1/2}= \sum_{i,j=1}^nq_i a^{1/2}q_j\in J^{1/2} $$ and hence $a\in J$.\ep
Notice that if $n=\infty$, then $\Psi(a)\in J$ does no longer imply that $a\in J$, as we see by considering non-compact positive matrices with  a compact main diagonal.

\begin{definition}\label{D5.3} Given $a\in M^+$, let
\begin{alignat*}{2}
a_-&:= (I-a)\chi_{a}(0, 1) \qquad &&\text{the defect  operator associated with $a$}\\
a_+&:= (a-I)\chi_{a}(1, \infty) &&\text{the excess  operator associated with $a$}.
\end{alignat*}
\end{definition}
Notice that  \be{e:5.5}a= a_+-a_-+R_a.\ee
When $M$ is a $\sigma$-finite factor, a characterization of strong sum of projections in $M$, i.e., of the sums of possibly infinitely many projections in $M$, with the convergence in the strong operator topology, was obtained in terms of the defect and excess operators in a previous work by three of the authors:

\bT{T5.4} \cite[Theorem 1]{KNZW*09}  Let $M$ be a $\sigma$-finite factor.
\item[(i)] Let  $M$ be of type I. Then $a$ is a strong sum of projections
 if and only if either $\tr(a_+)=\infty$ or \linebreak $\tr(a_-) \le \tr(a_+) < \infty$
 and $\tr(a_+) - \tr(a_-)\in \mathbb N\cup \{0\}$.
\item [(ii)]  Let $M$ be of type II and $a$ be diagonalizable (i.e., $a = \sum_{j=1}^\infty \gamma_n e_n$ for $\gamma_n\ge 0$ and $e_n\in M$ mutually orthogonal projections). Then $a$
is a strong sum of projections if and only if $\tau(a_+)\ge
\tau(a_-)$ where $\tau$ is a faithful semifinite normal trace on $M$. The condition is necessary even when $a$ is not
diagonalizable.
\item[(iii)] Let $M$  be of type III. Then  $a$ is a strong sum of projections if and only if either $||a||>1$ or $a$ is a projection.
\eT

We will present a necessary condition for $a$ to be a finite sum of projections in the case when $M$ is properly infinite. This condition is new even in the case of $B(H)$.

Recall that $J(M)$ is the closed two-sided ideal generated by the finite  $M$.  Denote by $F(M)$ the \textit{finite rank} ideal, namely the generally non-closed ideal of the elements of $M$ that have finite range projection. Then $F(M)$ and $J(M)$ share the same collection of projections, namely the collection of all the finite projections of $M$. 
\bT{T5.6}
Let $M$ be a  von Neumann algebra, $a\in M^+$, and $R_a$ be properly infinite. Assume that $a$ is a finite sum of projections in $M$. Then
\item [(i)] $\tau(a_+)\ge \tau(a_-)$ for any  (positive) trace $\tau$ on $M^+$ if $M$ is semifinite;
\item[(ii)] $c(R_{a_+})\ge c(R_{a_-})$;
\item [(iii)] if  $a_+ \in J(M)$, then $a_- \in J(M)$ and $\big( c(R_{a_+})- c(R_{a_-})\big )a_+\in F(M)$; 
\item [(iv)] if  $M$ is a direct sum of finitely many  factors, $a_+ \in J(M)$, and $a_-\ne 0$,  then   $a_+$ and $a_-$ generate the same two-sided (non-closed) principal ideal  of $M$. 
\eT
\bp
\item [(i)] By passing  if necessary  to $M_{R_a} $, assume without loss  of generality that $R_a=I$ and that $M$ is properly infinite. Then condition (iii) of Proposition \ref {P5.1} is satisfied and hence by condition (i) ibid there is a decomposition  of the identity $I=\sum_{j=1}^n q_j$ into mutually orthogonal projections $q_j\in M$  and an isometry $v$ such that $q_jvav^*q_j=q_j$ for all $j$. Equivalently, in the notations of Lemma \ref {L5.2}, 
$$\Psi(vav^*)=I.$$
Then by (\ref {e:5.5}) 
$$
\Psi(I)= I=\Psi(vav^*)=  \Psi(v(a_+-a_-+I) v^*)= \Psi(va_+ v^*)- \Psi(va_- v^*)+ \Psi(vv^*)
$$
and hence
\be{e:5.6}
 \Psi(va_+ v^*)= \Psi(va_- v^*)+ \Psi(I- vv^*).
 \ee
 Notice that the three operators $\Psi(va_+ v^*)$, $\Psi(va_- v^*)$, and $\Psi(I- vv^*)$ are positive.

 Thus if $\tau$  is trace    on $M^+$, we have $$\tau(a_+)=\tau(va_+v^*)= \tau\big(\Psi(va_+ v^*)  \big) \ge \tau\big(\Psi(va_- v^*)  \big)= \tau(va_-v^*)=\tau(a_-).$$

\item [(ii)] Since  $\Psi(va_+ v^*)\ge  \Psi(va_- v^*)$, and $c(R_{a_+})^\perp$ is central, we have
$$0= \Psi(vc(R_{a_+})^\perp a_+ v^*)=c(R_{a_+})^\perp\Psi(va_+ v^*)\ge c(R_{a_+})^\perp\Psi(va_- v^*)= \Psi(vc(R_{a_+})^\perp a_- v^*)\ge 0.$$
Thus $\Psi(vc(R_{a_+})^\perp a_- v^*)= 0.$ Since $vc(R_{a_+})^\perp a_-v^*\ge0$ and   $\Psi$ is faithful,  it follows that $vc(R_{a_+})^\perp a_-v^*=~0.$ Since  $v$ is an isometry, $c(R_{a_+})^\perp a_-=0$ and hence $c(R_{a_+})\ge c(R_{a_-})$.
\item [(iii)] Since $a_+\in J(M)$ it follows that $\Psi(va_+v^*)\in J(M)$. Since  $\Psi(va_-v^*)\le \Psi(va_+v^*)$ and  $\Psi(I-vv^*)\le \Psi(va_+v^*)$ and hence  it follows that $\Psi(va_-v^*)\in J(M)$ and $\Psi(I-vv^*)\in J(M)$. Then  by  Lemma \ref{L5.2}, 
$va_-v^*\in J(M)$ and $I-vv^*\in J(M)$.  Since $v$ is an isometry,  $a_-\in J(M).$ Furthermore, all projections in $J(M)$ are finite, hence 
\be{e:5.7}
I-vv^*\in F(M).
\ee 
By (ii), $c(R_{a_+})-c(R_{a_-})$ is a central projection and  
$$\big(c(R_{a_+})-c(R_{a_-})\big)\Psi(va_-v^*)= \Psi \Big(v\big(c(R_{a_+})-c(R_{a_-})\big) a_-v^*\Big)=0.
$$
Thus by (\ref {e:5.6}) and (\ref {e:5.7}) we have
$$
\Psi\Big(v\big(c(R_{a_+})-c(R_{a_-})\big)a_+v^*\Big)=  (c(R_{a_+})-c(R_{a_-}))\Psi(va_+v^*)= (c(R_{a_+})-c(R_{a_-}))\Psi(I-vv^*)\in F(M)
$$ 
whence we conclude by the same reasoning that $\big(c(R_{a_+})-c(R_{a_-})\big)a_+\in F(M)$.
\item [(iv)] The same argument as in (iii) shows that if $a_+$ belongs to some two-sided ideal $J$, then also $a_-\in J.$ Conversely, assume that  $a_-\in J$ for some nonzero ideal $J$.  If  $M$ is a factor, we have that $F(M)\subset J$ (see \cite[Theorem 6.8.3] {KrRj2}) and the same hold if $M$ is a direct sum of finitely many factors. Since by (\ref {e:5.7}), $I-vv*\in J$ and hence $\Psi(I-vv^*)\in J$, we have by (\ref {e:5.6})
  that
$
\Psi(va_+v^*)\in J$ and hence by Lemma \ref {L5.2} that $va_+v^*\in J$ and thus $a_+\in J$.
 \ep
 
If $M$ has infinite dimensional center, then (iv) may be false. This is due to the fact that there are ideals, even principal ideals generated by projections with central support the identity, that do do contain $F(M)$. The situation is illustrated by the following example.

\begin{example}\label{E5.7}
Let $M:=\bigoplus_1^\infty B(H)$ where $H$ is an infinite dimensional separable Hilbert space and let $I=\sum_{j=1}^\infty e_j$ be a decomposition of the identity into mutually orthogonal rank-one projections $e_j$. For every $n\ge 3$  set
$$ 
a_n:= \frac{1}{2}e_1+ \frac{3}{2}e_2+2\sum_{j=3}^ne_j+ \sum_{j=n+1}^\infty e_j.
$$
By Fillmore's theorem \ref {T1.1},  $\frac{1}{2}e_1+ \frac{3}{2}e_2$ is the sum of two projections, say $p_1+p_2$, and hence $a_n$ is also the sum of two projections $p_1+ \sum_{j=3}^n e_j$ and $p_2+ \sum_{j=3}^\infty e_j$. Thus $a:=\bigoplus_3^\infty a_n$ is also the sum of two projections.  Now, 
\ba
a_+&=\bigoplus_3^\infty \big(\frac{1}{2}e_2+\sum_{j=3}^ne_j\big)\\
a_-&=\bigoplus_3^\infty \frac{1}{2}e_1\qquad\text{and}\qquad c(a_-)=I.
\end{align*}
Thus 
$$a_+\in F(M)  = \{a=\bigoplus_1^\infty a_n\in M\mid a_n\in F(B(H)) ~~\forall n\}.$$
On the other hand, it is easy to see that the principal ideal generated by $a_+$ is $F(M)$ and the principal ideal generated by 
$a_-$ is 
$$
F_o(M):=  \{a=\bigoplus_1^\infty a_n\in F(M) \mid \sup_n\textnormal{rank}(a_n)< \infty\}
$$ 
and $F(M) \not\subset F_o(M)$.
Thus $a_+J(M)$, $a_-\ne 0$, but $a_+$ and $a_-$ do not generate the same principal ideal.
\end{example}

In the special case of $B(H)$ we have
 \bC{C5.8}
 Let $H$ be an infinite-dimensional separable Hilbert space, and assume that $a\in B(H)^+$ is a finite sum of projections and  that $a_+\in K(H)$. Then
 \item [(i)]  $a_-\in K(H)$.
  \item [(ii)] If $a_-=0$, then $a_+$ has finite rank.
   \item [(iii)] If $a_-\ne 0$, then $a_+$ and $a_-$ generate the same principal ideal in $B(H)$.
 \eC
\bC{C5.9}
Let $H$ be an infinite-dimensional separable Hilbert space  let $k_1 $ and $k_2$ be nonzero positive compact operators.
\item [(i)] $a:=I+k_1$ is a finite sum of projections if and only if $k_1$ has finite rank and integer trace.
\item[(ii)] If $b:=I+(k_1\oplus -k_2)$ is a finite sum of projections, then  $k_1$ and $k_2$ generate the same principal ideal in $B(H)$ and  either $\tr(k_1)=\infty$ or $\tr(k_1)-\tr(k_2)\in \mathbb N\cup\{0\}.$
\eC
 \bp
 \item [(i)] By definition, $a_-=0$ and $a_+=k_1$. Thus if $a$ is a finite sum of projections, by  Corollary \ref {C5.8} (ii) $k$ must have finite rank. By \cite[Theorem 1]{KNZW*09} (see  Theorem \ref {T5.4}) $\tr(k_1)=\tr (a_+)-\tr(a_-)\in \mathbb N\cup \{0\}$. Conversely, if $k_1$ has finite rank and integer trace, then $R_{k_1}+k_1$ is a finite sum of projections by Theorem \ref{T1.1}) and hence so is $a=  R_{k_1}^\perp +R_{k_1}+k_1$.
 \item [(ii)] By definition,  $b_+= k_1$ and $b_-= k_2$, hence the conclusion follows from Corollary \ref{C5.8} (iii) and \cite[Theorem 1]{KNZW*09} (see  Theorem \ref {T5.4}).
 \ep
 
 \begin{example}\label{E5.10}
Corollary \ref {C5.9} permits to construct simple examples of operators that are infinite sums of projections but not finite sums of projections. \\
In (i), choose $k_1$  with infinite rank and either infinite trace or integer trace. \\
In (ii) choose $k_1$ and $k_2$  that generate different principal ideals but that have the same trace (finite or infinite) and hence by Theorem \ref {T5.4} are strong sums of projections.  To construct such operators, we can choose  $k_1:=\diag(\xi_n)$ and $k_2:=\diag(\eta_n)$ with both sequences $\xi_n, \eta_n\downarrow 0$,  $\sum_1^\infty \xi_n= \sum_1^\infty \eta_n\le \infty$ and satisfying the $\Delta_{1/2}-$ condition, i.e., $\sup \frac{\xi_n}{\xi_{2n}}< \infty$, and $\sup \frac{\eta_n}{\eta_{2n}}< \infty$. Then $k_1$ and $k_2$ generate the same ideal if and only if $\alpha \xi_n\le \eta_n\le  \beta \xi_n$ for all $n $ and for some $0< \alpha\le \beta$. For instance, we can choose $\xi_n= \frac{1}{n}$ and $\eta_n= \frac{1}{n^{1/2}}$ or   $\xi_n= \frac{1}{2^n}$ and $\eta_n= \frac{6}{\pi^2n^2}$. For  background information on operator ideals in $B(H)$ we refer the reader to \cite  {DFWW}  (see also \cite {vKgW04-Density}).
\end{example}
\bR{R5.11}
There is some overlap between Corollaries \ref {C5.8}, \ref {C5.9}, and Example \ref {E5.10} and certain  results in \cite {WpCm09}. 
\item [(i)]
If  $a\in B(H)^+$ is a finite sum of projections and $\|a\|_e=1$, then $a_+\in K(H)$ and hence $a_-\in K(H)$ by Corollary  \ref {C5.8} (i). Thus $a - R_a = a_+-a_-\in K(H)$. This is \cite [Lemma 3.1]  {WpCm09}. 
\item [(ii)]
If $a$ is a finite sum of projections and $a_+-a_-\in K(H)$ has infinite rank, then by Corollary  \ref {C5.8} (ii) and (iii) both  $a_+$ and $a_-$ have infinite rank. This is the content of \cite [Theorem 3.4]{WpCm09}. 
\item [(iii)]
Corollary \ref {C5.9} (ii) provides a generalization of  \cite [Theorem 3.6] {WpCm09}. 
\item [(iv)]
Furthermore,  Example \ref {E5.10} is similar to, but  more general than  \cite [Example 3.9] {WpCm09}.
\eR

The problem of characterizing the positive operators of $B(H)$ that are finite sums of projections remains open.
\begin{question}\label{Q:12}
Is there a natural necessary and sufficient condition for the operator $b=I+(k_1\oplus -k_2)$ of Corollary \ref {C5.9} to be a finite sum of projections?
\end{question}
Notice that the technique used in \cite {KNZW*09} to decompose such an operator  $b$ into an infinite sum of rank one  projections, does not seem to provide a natural way to assemble those projections into a finite sums of projections. 

\section{Finite sums of projections in a type II$_1$ factor}
  In \cite {KNZW*09} (see Theorem  \ref {T5.4}  here)  we have proven that if $M$ is a  type II$_1$  factor $M$ with canonical trace $\tau$ and $a\in M$  is a strong sum of projections, then
\be{e:6.8}
\tau\big((a-I)\chi_a(1, \infty)\big) \ge \tau\big((I-a) \chi_a(0, 1)\big).
\ee
This condition is also sufficient in the cases when  $a$ is diagonalizable, that is $a=\sum_n\gamma_ne_n$ for some family of mutually orthogonal projections $e_n\in M$ and scalars $\gamma_n>0$. In the latter case, the  key step in the proof  \cite [Lemma 5.1]{KNZW*09} was the special case when 
\be{e:6.9} a= (1-\lambda)f+ (1+\mu)e\ee
 for a pair of mutually orthogonal  projections $e$ and $f$  in $M$ and scalars $0\le \lambda \le 1$ and $\mu>0$. For this case,  condition (\ref {e:6.8}) can be rewritten as 
\be{e:6.10}
\mu\tau(e) \ge  \lambda \tau(f)
\ee
and the proof in \cite [Lemma 5.1]{KNZW*09} was reduced to the case when equality holds in (\ref {e:6.10}).

Using a different approach we will show  that if $M$ is a type II$_1$ factor and $a$ is diagonalizable and strict inequality holds in (\ref{e:6.8}) then $a$ is a finite sum of projections.  We will start again with the special case (\ref {e:6.9}).

As in Lemma \ref {L4.2},  the strategy of the proof  is to decompose $a$ into a sum of two positive operators,  $b_o$ and the ``remainder" $a_1$ with the same form as $a$.     $b_o$ is  a positive combination of mutually orthogonal equivalent projections constructed so that $b_o$ satisfies Fillmore's necessary and sufficient condition (see Theorem \ref{T1.1}) for being a finite sum of projections. The construction is then  iterated and the  crux is to establish a bound on the number of projections needed in each block, so to satisfy the conditions of Lemma \ref {L4.1}.

   \bL{L6.1}
Let $M$ be  a type II$_1$ factor, $e$ and $f$ be mutually orthogonal projections, $0 \le  \lambda \le  1$, and $\mu>0$, and  $a: =  (1-\lambda)f+ (1+\mu)e$. If $\mu\tau(e) > \lambda \tau(f)$, then $a$ is a finite sum of projections in $M.$ Furthermore, there is an upper bound on the number of projections needed that depends only on $\mu$ if $f=0$ and on $\mu$ and $\frac{\mu\tau(e)}{\tau(f)}- \lambda$ if $f\ne 0$.
\eL

\bp
The cases when $\lambda=1$ (resp., $\lambda=0$) coincide with (resp.,  can be immediately reduced to) the case when $a= (1+\mu)e$. Hence assume that $0< \lambda < 1$.
Assume first that $\mu$ is not an integer. We will handle the case when $\mu$ is an integer at the end of the proof as a simple consequence of the construction in the non-integer case.
We construct recursively infinite sequences of projections $e_j$ and $f_j$ in $M$ and scalars $\lambda_j$, starting with  $e_{-1}:=e+f$,  $e_o: =e$, $f_o:=f$,  $\lambda_o:=\lambda$,  and setting   for all $j\ge 0$
\begin{align}
a_j:&=(1-\lambda_j)f_j+(1+\mu_j)e_j\label {e:6.11}\\
b_{j}:&=a_{j}-a_{j+1} \label {e:6.12}\\
\mu_{j+1}:&= \begin{cases}\mu_{j}- \lfloor \mu_{j} \rfloor \quad &f_{j}=0\\
\mu_{j} \quad &f_{j}\ne 0
 \end{cases} \qquad\text{starting with }\mu_o=\mu~~ \text{  ($\lfloor x \rfloor$ denotes the integer part of $x$)} \label{e:6.13}\\
 \delta_j:&= \frac{\mu_j\tau(e_j)}{\tau(f_j)}- \lambda_j ~\quad  f_j\ne0 \label {e:6.14}\\
\gamma_j:&=\min\{ \frac{\delta_j}{2}, \frac{1}{2}\}\qquad f_j\ne0 \label{e:6.15}
\end{align}
so that the
following conditions are satisfied for all $j\ge 0$
\begin{align}
&e_j\ne 0\label{e:6.16}\\
&e_jf_j= 0 \label{e:6.17}\\
&f_j+e_j\le e_{j-1}\label{e:6.18}\\
& 0< \lambda_j< 1\label{e:6.19}\\
&\frac{1}{2}\le \delta_j\le 2+ \mu \quad \text{for all } j\ge 1~ \text{for which  }f_j\ne0 \label {e:6.20}\\
& \mu_j\tau(e_j) > \lambda_j\tau(f_j)\label{e:6.21}\\
& \tau(e_{j+1})\le \big( 1-\frac{\mu_j-\lfloor\mu_j\rfloor}{2}\big)\tau(e_j) \label{e:6.22}\end{align}
\begin{align}
&R_{b_j} e_{j+1}=0 \quad\text{if } \mu_j= \mu_{j+1}\label{e:6.23}\\
& 0\ne b_j ~~ \text{is a sum of $N_j$ projections in $M$}  \label {e:6.24}\\
&N_j\le \begin{cases} 1+\mu& j=0,~ f=0\\ 
\frac{4(2+\mu)^2}{\mu}\max\{\delta_o,\frac{1}{\delta_o}\} & j=0,~ f\ne0\\
\frac{4(2+\mu)^3}{\mu}&j\ge 1
\end{cases}\label {e:6.25}
\end{align}
Once this construction is achieved,  we see from condition (\ref {e:6.13}) that $\mu_j\ne \mu_{j+1}$ for at most one index $j_o.$  Indeed, if $\mu>1$ and if $f_{j}= 0$ for some $j\ge0$, set $j_o$ to be the first such index. Then
\be{e:6.26}\mu_j=\begin{cases}\mu &j\le j_o\\\mu -\lfloor\mu\rfloor & j>j_o.\end{cases}\ee
If $\mu< 1$ or if $f_j\ne 0$ for all $j$, then $\mu_j=\mu$ for all $j$ and in this case, set $j_o=\infty$.

But then $$\mu_j-\lfloor\mu_j\rfloor =\mu-\lfloor\mu\rfloor\ne 0 \qquad \text{for all $j$,}$$ hence condition (\ref {e:6.22}) implies that $\tau(e_j)\to 0$ and thus $e_j\underset{s}\to 0$. From (\ref  {e:6.18}),  it follows that  also $f_j\underset{s}\to 0$ and since  the sequences $\{1-\lambda_j\} $ and $\{1+\mu_j\} $ are bounded by (\ref {e:6.19}) and (\ref{e:6.13}), it follows that $a_j\underset{s}\to 0$. By (\ref{e:6.12}), $a= \sum_{i=o}^ j b_i+a_{j+1}$ for all $j$  and hence $$a\underset{s}=\sum_{j=0}^\infty b_j.$$ Now notice that for all $j\ge 0$,
 \begin{alignat*}{2}R_{b_j}&\le R_{a_j} & & \text{(by (\ref  {e:6.12}))}\\
 &= f_j+e_j\hspace {3cm}&&\text{(by (\ref  {e:6.11}))}\\
 &\le e_{j-1}&&\text{(by (\ref  {e:6.18}).)} \end{alignat*}
Whenever $\mu_j= \mu_{j+1}$, i.e., $j\ne j_o$, then by (\ref {e:6.23}) we have $$R_{b_j}\perp e_{j+1}\ge R_{b_{j+2}}.$$
By (\ref {e:6.25}) and (\ref {e:6.26}), $\sup_j N_j=N< \infty$. At the end of the proof, we will find an explicit upper bound for $N$ including also the simple case when $\mu\in \mathbb N$.

 Thus the sequence $\{b_j\}_{j\ne j_o}$ satisfies the conditions of Lemma \ref {L4.1},  hence $\sum_{j\ne j_o} b_j$  is a sum of at most $2N$ projections and thus  $a$ is a sum of at most $3N$ projections.

\medskip
For the initial step of the construction, namely the construction of $e_1,f_1, \lambda_1$, we need to consider separately the simpler case when $f=0$ and the key case when $f\ne 0$.

\bigskip

\textbf{Case ($\mathbf {f=0}$)}  This is the case when $a= (1+\mu)e$.  (\ref{e:6.13}) prescribes that $\mu_1:=\mu_o- \lfloor \mu_o \rfloor$ and we set  $\lambda_1:=1-\mu_1$, thus satisfying (\ref{e:6.19}).
Since $M$ is of type II$_1$, we can choose  a projection $f_1\le e_o$ in $M$ with
$\tau(f_1)=\frac{\mu_1\tau(e_o)}{2}$ and set $e_1:=e_o-f_1$, thus satisfying (\ref{e:6.17}) and (\ref {e:6.18}).
Furthermore
$\tau(e_1)= (1-\frac{\mu_o-\lfloor\mu_o\rfloor}{2})\tau(e_o)$,
thus satisfying also (\ref{e:6.16}), and (\ref {e:6.22}).
Then
\begin{alignat*}{2}
\mu_1\tau&(e_1) - \lambda_1\tau(f_1) = \frac{\mu_1}{2}\tau(e_o)=\tau(f_1)>0 & \text{(thus satisfying (\ref{e:6.21}))} \\
a_1&=(1-\lambda_1)f_1+(1+\mu_1)e_1 &\text{(as prescribed by (\ref {e:6.11}))}\notag\\
&=(1+\mu_1)e_o-f_1\notag\\
b_o:&= a_o-a_1&\text{(as prescribed by (\ref {e:6.12}))}\notag\\
&= \lfloor \mu_o\rfloor e_o+ f_1& \text{(thus satisfying (\ref {e:6.24}) and (\ref {e:6.25}) with $N_o= \lfloor \mu_o\rfloor +1$ )}\notag\\
R_{b_o}&= \begin{cases} e_o&\mu_o>1\\
f_1&\mu_o<1.
\end{cases}\notag
\end{alignat*}
If $\mu_1=\mu_o$, i.e., $\mu_o<1$, then $ R_{b_o}=f_1\perp  e_1$, which satisfies (\ref {e:6.23}).
Thus conditions (\ref{e:6.16})--(\ref {e:6.22}) are satisfied.

Furthermore, $f_1\ne 0$ and
$\delta_1 =\frac{ \mu_1\tau(e_1) - \lambda_1\tau(f_1) }{\tau(f_1)}=  1 $, thus satisfying (\ref {e:6.20}).

\bigskip

\textbf{Case ($\mathbf {f\ne0}$)} .
Now we proceed with the first step of the construction in the key case when $f_o=f\ne0$. Define  nonnegative integers $k, n,$ and $ m$  and the positive real number $\alpha$ as follows:
\begin{alignat}{2}
k:&= \Big \lfloor \frac{2+\mu_o}{\gamma_o} \Big\rfloor &\text{i.e.,}\quad   \frac{2+\mu_o-\gamma_o}{\gamma_o}< k\le \frac{2+\mu_o}{\gamma_o}.\label{e:6.27}
\intertext{Then $k\ge 4$ since $\mu_o>0$ and $\gamma_o\le \frac{1}{2}$ and hence $ \frac{2+\mu_o}{\gamma_o}>4$.}
n:&=  \Big\lceil \frac{ k(\lambda_o + \delta_o -\gamma_o)}{\mu_o} \Big\rceil -1&\text{i.e.,} \quad \frac{ \lambda_o + \delta_o -\gamma_o}{\mu_o} \le \frac{n+1}{k} < \frac{ \lambda_o + \delta_o -\gamma_o}{\mu_o}+\frac{1}{k}. \label{e:6.28}
\intertext{Then $n\ge 0$ because $\frac{ k(\lambda_o + \delta_o -\gamma_o)}{\mu_o}>0$ since $\lambda_o>0$, $\delta_o -\gamma_o >0$, and $k>0$.}
m: &= \lfloor  (n+1)\mu_o-k \lambda_o  \rfloor+1\qquad\qquad &\text{i.e.,}\quad(n+1)\mu_o-k \lambda_o < m\le (n+1)\mu_o-k \lambda_o +1.\label{e:6.29}\\
\intertext{Then $m\ge 1$ because by (\ref{e:6.28}) we have $\frac{n+1}{k}\ge \frac{ \lambda_o + \delta_o -\gamma_o}{\mu_o}> \frac{ \lambda_o}{\mu_o}$ and hence $(n+1)\mu_o-k\lambda_o> 0$.}
\alpha:&= k \lambda_o -n \mu_o +m.\label{e:6.32}
\intertext{Then by (\ref {e:6.29})}
\mu_o& < \alpha \le 1+\mu_o.\label {e:6.33}
\intertext{As a consequence,}
(n+&1)\frac{\tau(f_o)}{k}< \tau(e_o).\label{e:6.34}
\end{alignat}
Indeed
\begin{alignat}{2}
 (n+1)\frac{\tau(f_o)}{k} &< \Big(  \frac{ \lambda_o + \delta_o -\gamma_o}{\mu_o}+\frac{\gamma_o}{2+\mu_o-\gamma_o}  \Big) \tau(f_o)\hspace{2.2cm}&(\text{by (\ref{e:6.28}) and (\ref{e:6.27})})\notag\\
&=\Big( \frac{ \lambda_o + \delta_o}{\mu_o} + \gamma_o\big(\frac{1}{2+\mu_o-\gamma_o} - \frac{1}{\mu_o}\big) \Big) \tau(f_o)\notag\\
&<\frac{ \lambda_o + \delta_o}{\mu_o}  \tau(f_o)&(\text{since $0< \gamma_o \le \frac{1}{2} <2$}) \notag\\
&= \tau(e_o) &(\text{by (\ref{e:6.14})}.)\notag
\end{alignat}
Use the fact that $M$ is of type II$_1$ to partition $f_o=\sum_{j=1}^kp_j$ into $k$ mutually orthogonal projections $p_j$, each with trace $\tau(p_j)= \frac{\tau(f_o)}{k}$  and  to find  $n+1$ mutually orthogonal equivalent projections $q_i \le e_o$, each with the same trace  $\tau(q_i)= \frac{\tau(f_o)}{k}$.  Now set

\begin{alignat}{2}
e_1&:=e_o-\sum_{i=1}^{n+1}q_i &\text{(thus by (\ref {e:6.34}) $\tau(e_1)>0$, satisfying (\ref{e:6.16}))}\notag\\
f_1&:= \begin{cases} q_{n+1} ~&\text{if }~\alpha \ne 1+\mu_o\\
0  ~&\text{if }~\alpha = 1+\mu_o
\end{cases}&\text{(thus satisfying (\ref{e:6.17}) and (\ref{e:6.18}))}\notag\\
\lambda_1&:= \begin{cases}\alpha-\mu_o ~&\text{if }~\alpha \ne 1+\mu_o\\ \text{arbitrary in }(0,1) ~&\text{if }~\alpha = 1+\mu_o
\end{cases}
&\text{(thus satisfying (\ref {e:6.19}))}\notag\\
\mu_1&= \mu_o &\text{(as prescribed by (\ref {e:6.13}))}\notag\end{alignat}
\begin{alignat}{2}
a_1&= (1-\lambda_1)f_1+(1+\mu_1)e_1 &\text{(as prescribed by (\ref{e:6.11}))}\notag\\
&=(1+\mu_o)e_o-(1+\mu_o)\sum_{i=1}^n q_i -\alpha q_{n+1}\notag\\
b_o&=a_o-a_1&\text{(as prescribed by (\ref{e:6.12}))}\notag\\
&=(1-\lambda_o)f_o+(1+\mu_o)\sum_{i=1}^n q_i +\alpha q_{n+1}\notag\\
&=\begin{cases} \sum_{j=1}^k(1-\lambda_o)p_j+ \sum_{i=1}^n(1+\mu_o)q_i + \alpha q_{n+1}& \text{if } n> 0\notag\\
\sum_{j=1}^k(1-\lambda_o)p_j + \alpha q_{n+1}& \text{if } n = 0.
\end{cases}\notag
\end{alignat}
Then
\begin{alignat}{2}
R_{b_o}&= \sum _{j=1}^kp_j+  \sum _{i=1}^{n+1}q_i\notag\\
&=f_o+e_o-e_1 &\text {(thus satisfying (\ref {e:6.23}))}\notag\\
\intertext{and}
\frac{\tau(e_1)}{\tau(e_o)}&= 1- \frac{(n+1)\tau(f_o)}{k\tau(e_o)}\qquad\qquad\qquad\qquad\qquad&(\text{by the definition of $e_1$})\notag \\
&= 1- \frac{(n+1)\mu_o}{k(\lambda_o+\delta_o)}&(\text{by (\ref {e:6.14})})\notag\\
&\le 1- \frac{\lambda_o+\delta_o-\gamma_o}{\lambda_o+\delta_o}  \qquad&(\text{by (\ref{e:6.28})})\notag \\
&< \frac{\gamma_o}{\delta_o}&(\text{since $\lambda_o>0$}) \notag \\
&\le \frac{1}{2}&(\text{by (\ref{e:6.15})})\notag\\
&< 1-\frac{\mu_o-\lfloor\mu_o\rfloor}{2}&(\text{thus satisfying (\ref {e:6.22}).)}\notag
\end{alignat}
Since all the $k+n+1$ mutually orthogonal projections in the decomposition of $b_o$ have the same trace, they are equivalent and we can view $b_o$ as belonging to a copy of $M_{k+n+1}(\mathbb C)$ embedded (not unitally) in $M$. Under this identification,  $b_o$ has   rank$(b_o)= k+n+1$ and by (\ref {e:6.32}),
$$
\tr(b_o)= k(1-\lambda_o)+ n(1+\mu_o)+\alpha = k+n+m\ge k+n+1= \text{rank}(b_o).
$$
Since $k$, $n$, and $m$ are integers, so is $\tr(b_o)$. Then by Theorem \ref{T1.1},   $b_o$ is the sum of $N_o:= k+n+m$ (equivalent) projections in $M$ and condition (\ref {e:6.24}) is satisfied. We claim that \be {e:6.35} N_o\le \frac{4(2+\mu_o)^2}{\mu_o}\max\{\delta_o,\frac{1}{\delta_o}\}
\ee
thus satisfying (\ref {e:6.25}).
Indeed
\begin{alignat*}{3}
N_o&=k=n+m\\
\le  k+ n +(n+1)\mu_o -k\lambda_o +1 \qquad &(\text{by (\ref{e:6.29})})\notag\\
&= (n+1)(1+\mu_o) + (1-\lambda_o)k\notag\\
&< \Big(  \frac{k(\lambda_o+\delta_o-\gamma_o)}{\mu_o}+1 \Big)(1+\mu_o)+ (1-\lambda_o)k\qquad &(\text{by (\ref{e:6.28}))}\notag\\
&= k\big(\frac{(\delta_o-\gamma_o)(1+\mu_o)}{\mu_o} +\frac{\lambda_o}{\mu_o}+1\big)+ 1+\mu_o\\
&<  k\big(\frac{(\delta_o-\gamma_o)(1+\mu_o)}{\mu_o} +\frac{1}{\mu_o}+1\big)+ 1+\mu_o &\text{since } \lambda_o< 1\\
&= (1+\mu_o)\Big(\frac{k} {\mu_o}(1+\delta_o-\gamma_o)+1\Big) \notag\\
&\le(1+\mu_o)\Big(\frac{(2+\mu_o)(1+\delta_o-\gamma_o)} {\mu_o\gamma_o}+1\Big)&(\text{by (\ref {e:6.27}))}\\
&= \begin{cases}(1+\mu_o)\Big(\frac{(2+\mu_o)(1+2\delta_o)} {\mu_o}+1\Big)&\delta_o\ge 1\\
(1+\mu_o)\Big(\frac{(2+\mu_o)(1+\frac{2}{\delta_o})} {\mu_o}+1\Big)&\delta_o\le 1
\end{cases} &\text{(by (\ref {e:6.15}))}\\
&= (1+\mu_o)\Big(\frac{(2+\mu_o)(1+2\max\{\delta_o,\frac{1}{\delta_o}\})} {\mu_o}+1\Big)\\
&= \frac{2(1+\mu_o)}{\mu_o}\big((2+\mu_o)\max\{\delta_o,\frac{1}{\delta_o}\} +1+\mu_o\big)\\
&\le \frac{4(2+\mu_o)^2}{\mu_o}\max\{\delta_o,\frac{1}{\delta_o}\}
\end{alignat*}
It may be interesting to notice but not essential for the remainder of the proof that a similar argument yields the inequality
$$N_o>\big(\frac{\delta_o}{2\mu_o}+1\big)(\frac{4}{\delta_o}-1\big)-1 $$
that  shows that $N_o$ must  indeed  be ``large" both when $\delta_o$  is ``large"  and when  $\delta_o$ is ``small". Thus to control the number of projections $N_j$ in the iterated construction, we need  upper and lower bounds on $\delta_j$ as given by (\ref {e:6.20}.)

 To obtain such bounds for $\delta_1$ (in the case when $f_1\ne0$, i.e., when $\alpha < 1+\mu_o$) and at the same time verify that  then the strict inequality (\ref {e:6.21}) holds, we start  from the following identity
\begin{alignat}{3}
 \delta_1 &= \frac{\mu_1\tau(e_1)} {\tau(f_1)}- \lambda_1 &(\text{by  (\ref {e:6.14})}) \notag\\
&=\frac{\mu_o\Big(\tau(e_o)- \frac{n+1}{k}\tau(f_o)\Big)}{ \frac{\tau(f_o)}{k}}-\alpha+\mu_o\qquad &(\text{by the definition of  $e_1$,$f_1$, $\lambda_1$ and $\mu_1$})\notag\\
&= \frac{\mu_o k \tau(e_o)}{\tau(f_o)}-  (n+1)\mu_o -(k\lambda_o-n\mu_o+m) + \mu_o\qquad &(\text{by (\ref{e:6.32})})\notag\\
&= k\delta_o -m &(\text{by (\ref {e:6.14}).})\notag\\
\intertext{Thus}
\delta_1&< k\delta_o -(n+1)\mu_o+k\lambda_o \qquad&(\text{by (\ref{e:6.29})})\notag\\
&\le k\gamma_o \qquad&(\text{by (\ref{e:6.28})})\notag\\
&\le 2+\mu_o &(\text{by (\ref{e:6.27})})\notag\\
\intertext{and}
\delta_1&\ge k\delta_o -(n+1)\mu_o+k\lambda_o -1\qquad&(\text{by (\ref{e:6.29})})\notag\\
&=k(\lambda_o+\delta_o)-(n+1)\mu_o-1\notag\\
&> k(\lambda_o+\delta_o) - k( \lambda_o+\delta_o-\gamma_o) -\mu_o-1 \qquad &(\text{by (\ref{e:6.28})})\notag\\
&=k\gamma_o-\mu_o-1\notag\\
&\ge 1-\gamma_o &(\text{by (\ref{e:6.27})})\notag\\
&\ge \frac{1}{2}&(\text{by (\ref{e:6.15}).})\notag
\end{alignat}
Thus (\ref {e:6.20}) is satisfied.
In particular, $\delta_1>0$ and hence condition (\ref{e:6.21}) is satisfied when $f_1\ne 0$. Notice that  (\ref{e:6.21})  is trivially satisfied when  $f_1=0$ since $e_1\ne 0$ by (\ref {e:6.16}).

\bigskip

Notice now that the remainder  $a_1= (1-\lambda_1) f_1+ (1+\mu_1)e_1$  has precisely the same form as the original element $a_o$, but while $\delta_o$ (in the case when $f_o\ne 0$) was arbitrary, $\delta_1$ (in the case when $f_1\ne 0$)  is bounded above and below by (\ref {e:6.20}) and hence $\max\{\delta_1,\frac{1}{\delta_1}\}\le 2+\mu_1$. Thus we can apply to $a_1$ the construction of Case ($f=0$) if $f_1=0$ or Case ($f\ne0$) if $f_1\ne 0$ and obtain a better  estimate for the upper bound of $N_2$. Indeed,  from (\ref {e:6.35}) applied to this second step
$$N_2\le \begin{cases} 1+ \mu_1 &f_1=0\\
\frac{4(2+\mu_1)^2}{\mu_1}\max\{\delta_1,\frac{1}{\delta_1}\}&f_1\ne0
\end{cases} \quad \le \quad \frac{4(2+\mu_1)^3}{\mu_1} \le \max\Big\{\frac{4(2+\mu)^3}{\mu}, \frac{4(2+\mu-\lfloor\mu\rfloor)^3}{\mu-\lfloor\mu\rfloor}\Big\}
$$
thus satisfying (\ref {e:6.25}) for $j\ge 1.$
The recurrence then proceeds exactly as in the case $j=1$, thus concluding the proof for the case when $\mu\not \in \mathbb N$.

Finally, assume that $\mu\in \mathbb N$. If $f\ne 0$ we can apply the construction of Case ($f\ne 0$) without any changes. And if then $f_1\ne 0$ we can apply the construction again. There are two possibilities.

Either $f_{ j}$ never vanishes, i.e.,  the process continues indefinitely applying only the construction in Case ($f\ne 0$), in which case $\mu=\mu_j$ for all $j$ and the same conclusion and the same bounds for $N_j$ as in the first part of the proof hold,  or there is some $j$ for which $f_{ j}=0.$ Let $j_o$ be the first such index. 
In that case, $$a_{j_o}= (1+\mu_{j_o})e_{j_o}= (1+ \mu) e_{j_o}$$ is already the sum of $1+\mu$ projections  and the ``remainder"  $a_{j_o+1}$ vanishes. Thus the process terminates at step $j_o$. But for all $0\le j \le j_o$, $\mu_j=\mu$, the  bounds for $N_j$ (see (\ref {e:6.25}) and the relations between the range projections of the elements $b_j$ and the projections $e_j$ that we have  obtained in the first part of the proof hold, and hence, so does the conclusion.

Now to summarize, we have in all cases (i.e., whether $\mu$ is an integer or not),
\ba N_o&\le \begin{cases} 1+ \mu \hspace{3.4cm}&f=0\\
\frac{4(2+\mu)^2}{\mu}\max\{\delta_0,\frac{1}{\delta_0}\}&f\ne 0
\end{cases}\\
\sup_{j\ge 1} N_j&\le \begin{cases}\frac{4(2+\mu)^3}{\mu }&j_o=\infty\\ 
\max\Big\{\frac {4(2+\mu)^3}{\mu}, \frac{4(2+ \mu- \lfloor\mu\rfloor)^3}{ \mu- \lfloor\mu\rfloor }  \Big\} &j_o<\infty. \end{cases}
\end{align*}
But then, the bound $N=\ max\{N_o, \sup_{j\ge 1} N_j\} $ on the number of needed projections for each block, and hence the bound $3N$ on the number of projections needed to decompose $a$ is seen to indeed depend only on $\mu$ if $f\ne 0$, and on $\mu$ and $\delta_o= \frac{\mu \tau(e )}{\tau(f )}- \lambda $ when $f\ne 0$.
\ep
\bR{R6.2}
\item[(i)] If $a= (1+\frac{m}{n})e$ with $m, n\in \mathbb N$, then instead of proceeding with the Case (f=0) construction etc, we can decompose directly  $e$ into $n$ mutually orthogonal equivalent projections. Then by identifying $a$ with an element of $M_n(\mathbb C)$, we see that rank$(a)=n$ and $\tr(a)= n+m$, hence $a$ is the sum of $n+m$ (equivalent) projections in $M$. In other words, if $\mu$ is rational we can terminate the construction in Lemma \ref {L6.1} at any step where $f_j=0$. 
\item[(ii)] If $a=(1-\lambda)f+ (1+\mu)e$ and $\frac{\mu}{\lambda}= \frac{\tau(f)}{\tau(e)}= \frac{m}{n}$ for some $m,n\in \mathbb N$, then by decomposing $e$ into $m$ mutually orthogonal equivalent projections  and $f$ into $n$ mutually orthogonal equivalent  the same reasoning as in (i) shows that $a$ is a sum of $n+m$ equivalent projections. This observation shows that the strict inequality $\tau(a_+) > \tau(a_-)$ is not necessary for $a$ to be a finite sum of projections.
\eR


\bigskip

Before we consider the general case of a positive diagonalizable operator with $\tau(a_+) > \tau(a_-)$, we need  the following elementary lemma.
\bL {L6.4} Assume that $\sum_{i=1}^m\xi_i=\sum_{j=1}^n\eta_j$ for some $m,n\in \mathbb N$ and $\xi_i, \eta_j>0$ for $1\le i\le m,~1\le j\le n$. Then for all $1\le i\le m,~1\le j\le n$ there are decomposition $\xi_i=\sum_{k=1}^{m_i}\xi_{i,k}$, $\eta_j=\sum_{h=1}^{n_j}\eta_{j,h}$ where $\xi_{i,k}>0$, $\eta_{j,h}>0$ and such that $\{\xi_{i,k}\}_{1\le i \le m,~1\le k\le m(i)}= \{\eta_{j,h}\}_{1\le j \le n,~1\le h \le n(j)}.$
\eL
\bp
\item The proof is on induction on the number $m+n$ of summands. Assume without loss of generality that $m\le n$. To simplify notations, assume also that the finite sequences $\xi_i$ and $\eta_j$ are already in monotone non-increasing order. Then $m\xi_1\ge n\eta_n$. Thus if $\eta_n\ge \xi_1$ then $m=n$, $\xi_i=\eta_j$ for all $i$ and $j$ and there is nothing to prove. Thus assume that $\eta_n< \xi_1$ and let $\xi_1':=  \xi_1-\eta_n$. Then $\xi_1'+\sum_{i=2}^m\xi_i= \sum_{j=1}^{n-1}\eta_j$ satisfies the same hypotheses but with a number $m+n-1$ of summands. Routine arguments conclude the proof.
\ep
\bT{T6.5}
Let $M$ be a  type II$_1$ factor and let $a\in M$ be a positive  diagonalizable operator. Then a sufficient condition for  $a$ to be a finite sum of projections in $M$ is  that  $\tau((a-I)\chi_a(1, \infty)) > \tau((I-a) \chi_a(0, 1)).$
\eT
\bp
Assume without loss of generality that  $1$ is not a eigenvector of $a$ and hence  that
$$a= \sum_{j=1}^{N_1}(1-\lambda_j)f_j + \sum_{i=1}^{N_2}(1+\mu_i)e_i,
$$
where $N_1 \in \mathbb N\cup \{\infty\}\cup \{0\}$ (again, adopting the convention that a sum like $ \sum_{j=1}^{0}$ is  dropped), $N_2 \in \mathbb N\cup \{\infty\}$,  $\mu_i>0$, $0< \lambda_i< 1$ for all $i$ and $j$, and $\{e_i,f_j\}$ are mutually orthogonal nonzero projections in $M$ whose sum is the range projection $R_a$ of $a$. Then
$$a_+=(I-a)\chi_a(0, 1)=  \sum_{j=1}^{N_1}\lambda_jf_j\qquad\text{and}\qquad a_-=(a-I)\chi_a(1, \infty)= \sum_{i=1}^{N_2}\mu_ie_i.$$ Set
$$ \gamma:= \tau(a_+) - \tau(a_-)= \sum_{i=1}^{N_2}\mu_i\tau(e_i)- \sum_{j=1}^{N_1}\lambda_j \tau(f_j).
$$
Thus $\gamma >0$ by hypothesis. Assume that $N_2=N_1=\infty$ (the other cases are simpler and are left to the reader.)   Set $h=  \lceil 14 \|a\|\rceil -1 $. By hypothesis, $\|a\|>1$ and hence $h\ge 14$. Choose $N\in \mathbb N$ so that
$$
\sum_{i=N+1}^\infty\tau(e_i) \le \min\big\{ \frac{\gamma}{3(\|a\|-1)}, \frac{\gamma}{3h}\big\}.
$$
Let $e_\infty:= \sum_{i=N+1}^\infty e_i$ and $a_o: = \sum_{i=N+1}^\infty (1+\mu_i) e_i$. Then $a_o\ge e_\infty =R_{a_o}$ and we can apply  Corollary \ref {C2.4} with $\nu=1$, and by (\ref{e:2.1})  with $N_o=12$, $V_o=14$ . Thus $a$ can be expressed as a positive linear combination of not more than
$$N_o+\lceil V_o \big( \frac{ \|a_o\|}{\nu}-1\big)\rceil +1\le 12 +\lceil 14 \|a\| -14\rceil +1 =h$$ projections. Thus  $a_o= \sum_{k=1}^n\alpha_k p_k$ for some collection of  $n\le h$ nonzero projections $p_k\in M$ and positive scalars $\alpha_k$. Necessarily, $\lfloor \alpha_k\rfloor\le \alpha_k\le \|a_o\|\le \|a\|$.  By decomposing
$$a_o=  \sum_{k=1}^n\lfloor \alpha_k\rfloor p_k+  \sum_{k=1}^n\big(\alpha_k - \lfloor \alpha_k\rfloor\big)p_k$$
and noticing that the first summand is the sum of at most $n\|a\|$ projections, to simplify notations,  assume that $0< \alpha_k< 1$ for each $1\le k\le n$ . Set  $\beta_k=1-\alpha_k$. Thus
\be{e:6.36}a= \sum_{i=1}^{N}(1+\mu_i)e_i+ \sum_{j=1}^{\infty}(1-\lambda_j)f_j + \sum_{k=1}^{n}(1-\beta_k)p_k.\ee
Notice that the projections $p_k$ are not mutually orthogonal, but since $p_k\le e_\infty$ for all $k$, each $p_k$ is  orthogonal to all the projections $e_i$ for $1\le i\le N$ (and also to all the projections $f_j$, although we do not need the latter fact.)
Furthermore,
$$
\sum_{k=1}^{n}\beta_k\tau (p_k) < \sum_{k=1}^{n}\tau (p_k)\le n\tau(e_\infty)\le \frac{\gamma}{3}.
$$
We also have
$$
\sum_{i=N+1}^\infty\mu_i\tau(e_i)\le \sum_{i=N+1}^\infty(\|a\|-1)\tau(e_i) = (\|a\|-1)\tau(e_\infty) \le \frac{\gamma}{3}.
$$
But then,
$$
\sum_{i=1}^N \mu_i\tau(e_i) - \sum_{j=1}^\infty \lambda_j\tau(f_j)- \sum_{k=1}^{n}\beta_k\tau (p_k)=\gamma- \sum_{k=1}^{n}\beta_k\tau (p_k) - \sum_{i=N+1}^\infty\mu_i\tau(e_i)\ge \frac{\gamma}{3}>0.
$$
Set
\ba \rho&:= \sum_{i=1}^N \mu_i\tau(e_i) -  \sum_{j=1}^\infty \lambda_j\tau(f_j)-\sum_{k=1}^{n}\beta_k\tau (p_k)\\
\sigma&: = \sum_{j=1}^\infty \tau(f_j)+ \sum_{k=1}^{n}\tau (p_k).
\end{align*}
Then $\rho, \sigma >0$ and 
\begin{align}
\sum_{i=1}^N \mu_i\tau(e_i)&= \sum_{j=1}^\infty\lambda_j\tau(f_j)+\sum_{k=1}^{n}\beta_k\tau (p_k)+\rho\notag\\
&=\sum_{j=1}^\infty\lambda_j\tau(f_j)+\sum_{k=1}^{n}\beta_k\tau (p_k) +\frac{\rho }{\sigma}\Big( \sum_{j=1}^\infty\tau(f_j)  +\sum_{k=1}^{n}\tau (p_k) \Big)\notag\\
&=
 \sum_{j=1}^\infty\big ( \lambda_j+ \frac{\rho }{\sigma}\big)\tau(f_j)+\sum_{k=1}^{n}\big(\beta_k+ \frac{\rho}{\sigma}\big)\tau (p_k)\label{e:6.37}
\end{align}
Choose $m\in \mathbb N$ for which $  \sum_{j=m+1}^\infty\big ( \lambda_j+ \frac{\rho }{\sigma}\big)\tau(f_j)
\le \mu_1\tau(e_1)$. Since $M$ is of type II$_1$,  we can find mutually for each orthogonal projections $e_{1,j}\le e_1$ for all $j\ge m+1$ such that 
\be{e:6.38}
\tau(e_{1,j})= \frac{1}{\mu_1}\big ( \lambda_j+ \frac{\rho }{\sigma}\big)\tau(f_j).
\ee
Set for $j\ge m+1$
$$ a(j):=(1+\mu_1)e_{1, j}+(1-\lambda_j)f_j\qquad\text{and}\qquad
a' :=  \sum_ {j=m+1}^\infty a(j).
$$
Since by (\ref{e:6.38}) $$ \frac{\mu_1\tau(e_{1,j})}{\tau(f_j)}- \lambda_j= \frac{\rho}{\sigma}>0, $$
it follows that $a(j)$ satisfies the conditions of Lemma \ref {L6.1} and hence is the sum of finitely many projections. Moreover, for each $j$, there is an upper bound $N(j)$ on the number of projections needed that depends only on  $\mu_1$ and $\frac{\mu_1\tau(e_{1,j})}{\tau(f_j)}- \lambda_j$, both of which are constant for all $j$. Thus $N(j)$ is also constant, say $N(j)\equiv N'$.  As a consequence, $a'$ is also a sum of finitely many projections.

Now set \begin{align}
a''&:=a-a'\notag\\
e_1''&:= e_1- \sum_{j=m+1}^\infty e_{1,j}.\label{e:6.39}
\end{align}
Then 
\be{e:6.40}
a''=(1+\mu_1) e_1''+ \sum_{i=2}^N(1+\mu_i)e_i+\sum_{j=1}^m(1-\lambda_j)f_j+ \sum_{k=1}^{n}(1-\beta_k)p_k.\ee
Now we have the identity
\begin{alignat*}{2}
\mu_1\tau( e_1'')+ \sum_{i=2}^N \mu_i\tau(e_i) &=\mu_1\tau(e_1)- \mu_1 \sum_{j=m+1}^\infty\tau(e_{1,j})+\sum_{i=2}^N \mu_i\tau(e_i)\qquad &&(\text {by (\ref{e:6.39}))}\\
&=\sum_{i=1}^N \mu_i\tau(e_i) - \sum_{j=m+1}^\infty \big ( \lambda_j+ \frac{\rho }{\sigma}\big)\tau(f_j)\qquad &&(\text {by (\ref{e:6.38}))}\\
&=\sum_{j=1}^m \big ( \lambda_j+ \frac{\rho }{\sigma}\big)\tau(f_j)+ \sum_{k=1}^{n}\big(\beta_k+ \frac{\rho}{\sigma}\big)\tau (p_k) &&(\text {by (\ref{e:6.37}).)}\\
\end{alignat*}
Now apply Lemma \ref {L6.4} to this identity, that is, decompose each of the summands so to ``match" the subsummands and notice that 
every decomposition of the (scalar) summands in this identity leads to a corresponding decomposition in the projections. For instance,  if for some $i$ we have the decomposition $\mu_i\tau(e_i)= \sum_{k=1}^{m(i)} \xi_{i,k}$, then, using the fact that $M$ is a type II$_1$ factor we can correspondingly decompose $e_i$ into the sum $e_i=\sum_{k=1}^{m(i)} e_{i,k}$ of mutually orthogonal projections with traces $\tau( e_{i,k})= \frac{ \xi_{i,k}}{\mu_i}$. Similarly, decompose correspondingly all of the other projections $f_j$ and $p_k$.
To simplify notations, list as $\{f'_i\}_1^{m'}$ (resp., $\{p'_i\}_{m'+1}^{n'}$) all the projections obtained from the decomposition of the projections $\{f_j\}_1^m$ (resp., $\{p_k\}_{1}^{n}$) and therefore list as $\{e'_i\}_1^{n'}$ the projections obtained from the decomposition of the projections $e_1''$ and $\{e_j\}_1^N$ and label accordingly the corresponding coefficients. Explicitly, 
\be{e:6.41}
\mu'_i\tau(e'_i)=
\begin{cases}
\big(\lambda'_i+ \frac{\rho}{\sigma}\big) \tau(f'_i)   &1\le i\le m'\\
\big(\beta'_{i}+ \frac{\rho}{\sigma}\big)\tau(p'_{i})& m'+1\le i\le n'
\end{cases} 
\ee
Using this decomposition, we rewrite (\ref{e:6.40}) as 
\ba a''&= \sum_{i=1}^{n'}(1+\mu'_i)e'_i+ \sum_{i=1}^{m'}(1-\lambda'_i)f_i'+  \sum_{i=m'}^{n'}(1-\beta'_i)p'_i \\
&= \sum_{i=1}^{m'}\Big( (1+\mu'_i)e'_i+ (1-\lambda'_i)f'_i\Big)+  \sum_{i=m'+1}^{n'}\Big( (1+\mu'_i)e'_i+ (1-\beta'_{i})p'_{i}\Big)
\end{align*}
Thus $a''$ is now decomposed into the sum of  the $n'$ positive  operators 
$$\begin{cases}(1+\mu'_i)e'_i+(1-\lambda'_i)f_i'& 1\le i\le m'\\
(1+\mu'_i)e'_i+(1-\beta_{i})p'_{i'}& m'+1\le i\le n'.
\end{cases}$$

Recall that while the projections $p_k$ may not be mutually orthogonal, they are orthogonal to each of the projections $e_i$ and hence so are each of the projections $p_i'$ into which they have been decomposed are orthogonal to each of the projections $e_j'$. In particular $e'_i$ is orthogonal to $p_i'$ for every $m'+1\le i\le n'$ and of course, $e'_i$ is orthogonal to $f_i'$ for every $1\le i\le m'$. By
(\ref {e:6.41}) we have 
$$
\mu'_i\tau(e'_i)>
\begin{cases}
\lambda'_i\tau(f'_i) &1\le i\le m'\\
 \beta'_{i}\tau(p'_{i})& m'+1\le i\le n'.
\end{cases} 
$$
Thus by Lemma \ref {L6.1}, each of these positive operators
is a finite sum of projections,  and hence so is $a''$, This completes the proof.  \ep

An obvious consequence of this theorem is
\bC{C6.6}
If $M$ is a type II$_\infty$ factor, $a\in M^+$ is diagonalizable, $R_a$ is finite, and $\tau(a_+)> \tau(a_-)$, then $a$ is a finite sum of projections.
\eC
\bR{R6.7}
The condition that $\tau(a_+)>\tau(a_-)$ alone, (i.e., without requiring also that $R_a$ be finite,)  is not sufficient for $a$ to be a finite sum of projections (whether $a$ is diagonalizable or not).  For instance, if  $a_+ \in J(M)$, $\tau(a_+)=\infty $ (necessarily $R_a$ is infinite), but  $\tau(a_-)<\infty $, then  $a_-$ belongs to the trace-class ideal while and $a_+$ does not, so by Theorem \ref {T5.6}, $a$ is not a finite sum of projections.
\eR
We have used the condition that $a$ is diagonalizable in order to reduce the problem to the simpler case of Lemma \ref {L6.1}.
However,  this condition is clearly not necessary for $a$ to be a finite sum of projections. For instance, if $0\le a\le 2I$ and there is a unitary $ u$ for which $uau^*=2I-a$, then $a$, whether diagonalizable or not,  is sum of two projections by  \cite [Corollary of Theorem 2]{Fp69} (see \cite [Proposition 2.10] {KNZW*09} for the von Neumann algebra case.)
 \begin{question}\label{Q6.8}
Can the condition that $a$ is diagonalizable be removed?
\end{question}

\end{document}